\documentclass[letterpaper]{article}
\usepackage[intlimits]{amsmath}
\usepackage{amsfonts, amssymb, amsthm}
\usepackage{epsfig}
\usepackage{color}
\usepackage{pstricks}
\newlength{\hchng}
\newlength{\vchng}
\newtheorem{thm}{Theorem}[section]
\newtheorem{prop}[thm]{Proposition}
\newtheorem{cor}[thm]{Corollary}
\newtheorem{lemma}[thm]{Lemma}

\newtheorem{preremark}[thm]{Remark}
\newenvironment{remark}{\begin{preremark}\rm}{\medskip \end{preremark}}
\numberwithin{equation}{section}

\newcommand{\norm}[1]{\left\Vert#1\right\Vert}
\newcommand{\abs}[1]{\left\vert#1\right\vert}
\newcommand{\set}[1]{\left\{#1\right\}}
\newcommand{\R}{\mathbb R}
\newcommand{\eps}{\varepsilon}

\newcommand{\grad} {\nabla}
\newcommand{\lap} {\triangle}
\newcommand{\bdary} {\partial}
\newcommand{\dx} {\; \mathrm{d} x}
\newcommand{\dd} {\; \mathrm{d}}

\DeclareMathOperator*{\osc}{osc}

\DeclareMathOperator{\dist}{dist}
\DeclareMathOperator{\dv}{div}

\newcommand{\fl}{(-\lap)^s}

\newcommand{\fld}{\lim_{y \to 0^+} y^a \partial_y}
\newcommand{\Dy}[1]{\dv \left(y^a \grad #1\right)}
\newcommand{\Dey}[1]{L_a #1}

\newcommand{\yed} {\; |y|^a \mathrm{d}}

\def\XXint#1#2#3{\quad {\setbox0=\hbox{$#1{#2#3}{\int}$}
       \vcenter{\hbox{$#2#3$}}\kern-0.5\wd0}}

\title{Regularity estimates for the solution and the free boundary to the obstacle problem for the fractional Laplacian}
\author{Luis Caffarelli, Sandro Salsa and Luis Silvestre}
\begin{document}
\maketitle

\begin{abstract} We use a characterization of the fractional Laplacian as a Dirichlet to Neumann operator for an appropriate differential equation to study its obstacle problem. We write an equivalent characterization as a thin obstacle problem. In this way we are able to apply local type arguments to obtain sharp regularity estimates for the solution and study the regularity of the free boundary.
\end{abstract}

\section{Introduction}
Constrained variational problems with fractional diffusion appear in the study of the quasi-geostrophic flow model (\cite{Co} and \cite{CV}), anomalous diffusion \cite{BG} and American options with jump processes \cite{CT} (See also \cite{S}).

In this paper we will study the \emph{fractional} obstacle problem. It appears in several contexts and it can be stated in different ways. Given a smooth function $\varphi : \R^n \to \R$, that we can assume it decays rapidly at infinity, we look for a function $u$ satisfying
\begin{equation} \label{eq:flp}
\begin{aligned}
u(x) &\geq \varphi && \text{for } x \in \R^n,\\
\fl u(x) &= 0 && \text{for } u(x)>\varphi, \\
\fl u (x)&\geq 0 && \text{for } x \in \R^n.
\end{aligned}
\end{equation}

This system of inequalities may also be written as a variational problem in $\dot H^s$ (The Hilbert space spanned by $C_0^\infty$ functions with the norm $\norm{v}_{\dot H^s} = \norm{ \hat v (\xi) |\xi|^s }_{L^2}$). The function $u$ minimizes its norm in $\dot H^s$ among all functions $v$ satisfying $v \geq \varphi$.

From a potential theoretic point of view, $u$ can also be characterized as the smallest supersolution of $(-\lap)^s v \geq 0$, among those functions $v \geq \varphi$. This makes sense since the minimum of supersolutions is also a supersolution.

Finally, as a Hamilton-Jacobi equation we can describe $u$ by the property
\[ \min( (-\lap)^s u, u-\varphi) = 0 \]

Each of the previous descriptions is motivated by a different application and different interpretations suggest how to treat different issues as we study them.

In this paper we develop two aspects of the problem: optimal regularity of the solution and regularity of the free boundary. The existing quasi-optimal regularity of the solution $u$ was proven by one of the authors in \cite{S}. There it is shown that $u \in C^{1,\alpha}$ for every $\alpha \in (0,s)$ using methods mainly from potential analysis.

In the particular case $s=1/2$ it is easy to see that the operator $(-\lap)^{1/2}$ coincides with the Dirichlet to Neumann operator in the upper half space of $\R^{n+1}$. More precisely, given $u(x)$ defined in $\R^n$,extend it to $u^*(x,y)$ in $\R^{n+1}$ by convolving with the classical Poisson kernel. Then $(-\lap)^{1/2} u(x) = -u_y^* (x,0)$.

In \cite{CS}, we obtain the same interpretation for any fractional power $(-\lap)^s u$ ($0<s<1$) as the Dirichlet to Neumann operator of an appropriate extension $u^*(x,y)$. In this paper we use this characterization to obtain the sharp optimal regularity $u \in C^{1,s}$. We also obtain a regularity result for the free boundary away from singular points.

We obtain the equivalent problem to \eqref{eq:flp} extending $u$ to the upper half space $u : \R^n \times [0,\infty) \to \R$.
\begin{align}
u(x,0) &\geq \varphi(x) && \text{for } x \in \R^n \label{eq:ugp} \\
\Dy{u(x,y)} &= 0 && \text{for } y>0 \label{eq:usu} \\
\fld u(x,y) &= 0 && \text{for } u(x,0)>\varphi(x) \label{eq:uzn}\\
\fld u(x,y) &\leq 0 && \text{for } x \in \R^n \label{eq:uln}
\end{align}
where $ s = (1-a)/2$.

For $y>0$, $u(x,y)$ is smooth, thus the equation \eqref{eq:usu} is understood in the classical sense. The equations at the boundary \eqref{eq:uzn} and \eqref{eq:uln} should be understood in the weak sense as it is explained in \cite{CS}. Also, from \cite{S}, $u(x,0) \in C^{1,\alpha}$ for any $\alpha<s$, in particular for some $\alpha$ such that $2s < 1+\alpha < 1+s$, then $\fld u(x,y)$ can be understood in the classical sense too.

From \cite{S} we know that $u(x,0) \in C^{1,\alpha}$ for every $\alpha<s$, which easily implies the same for $u(x,y_0)$ uniformly with respect to any fixed $y_0$. We also know from \cite{S} that for $C = \sup \abs{D^2 \varphi}$, $\partial_{ee} u \geq -C$ for any unit vector $e$ ($u$ is semiconvex in the $x$ variable).

The function $u$ that solves the problem above can be extended to the whole space by symmetrization $u(x,y) = u(x,-y)$. In \cite{CS} it was shown that the condition \eqref{eq:uzn} is equivalent to the symmetric extension of $u$ being a solution of $\dv(|y|^a \grad u(X) ) = 0$ across the boundary $\{y=0\}$ on the part where $u(x,y)>\varphi(x)$. On the other hand, the condition \eqref{eq:uln} is equivalent to $\dv(|y|^a \grad u(X) ) \leq 0$ as a distribution. Let us call $L_a u = \dv( |y|^a \grad u)$. The setting for the symmetric extension translates as
\begin{align}
u(x,0) &\geq \varphi(x) && \text{for } x \in \R^n \label{eq:symrn1}\\
u(x,y) &= u(x,-y) \\
\Dey{u(X)} &= 0 && \text{for } X \in \R^{n+1} \setminus \{(x,0) : u(x,0) = \varphi(x) \} \\
\Dey{u(X)} &\leq 0 && \text{for } X \in \R^{n+1} \text{ in the distributional sense} \label{eq:symrn4}
\end{align}

Notice that $L_a u$ is a singular measure supported on the set $\Lambda := \{ (x,0) : u(x,0) = \varphi(x) \}$. By the continuity of $u$, we have that $(u(X) - \varphi(x)) L_a u(X) = 0$ in $\R^{n+1}$.

In the case $s=1/2$ (i.e. $a=0$), $L_a u = \lap u$, and we are on the situation of the classical thin obstacle problem. Optimal regularity estimates, as well as free boundary regularity results have been obtained recently for this problem in \cite{AC}, \cite{ACS}. However, the results in this paper are new even to that case because we can consider nonzero obstacles $\varphi$.

One of the great advantages of stating the problem with a PDE (\ref{eq:ugp}-\ref{eq:uln}) instead of a nonlocal equation \eqref{eq:flp} is that we can localize the problem. We will consider a local version of (\ref{eq:ugp}-\ref{eq:uln}). We write $X = (x,y)$, with $X \in \R^{n+1}$, $x\in \R^n$, $y \in \R$, and
\[ \begin{aligned}
B_r &= \{X: X \in \R^{n+1}, |X| < r \} \\
%B_r^+ &= \{(x,y): x \in \R^n, y \in [0,\infty), |x|^2+|y|^2 < r^2 \} \\
B_r^* &= \{x \in \R^n: |x| < r \} \\
S_r   &= \{X: X \in \R^{n+1}, |X| = r \}\\
%S_r^+   &= \{(x,y): x \in \R^n, y \in [0,\infty), |x|^2+|y|^2 = r^2 \}
\end{aligned}\]
% S_r may become S_r^+ to denote that it is only half sphere.

Given $\varphi :B_1^* \to \R$, we consider a function $u : B_1 \to \R$ satisfying the following equation
\begin{align}
u(x,0) &\geq \varphi(x) && \text{for } x \in B_1^* \label{eq:symugp}\\
u(x,y) &= u(x,-y) \label{eq:symsym}\\
\Dey{u(X)} &= 0 && \text{for } X \in B_1 \setminus \{(x,0) : u(x,0) = \varphi(x) \} \label{eq:symeq}\\
\Dey{u(X)} &\leq 0 && \text{for } X \in B_1 \text{ in the distributional sense} \label{eq:symineq}
\end{align}

Compared to (\ref{eq:symrn1}-\ref{eq:symrn4}), we are forgeting everthing that happens outside $B_1$. We can do that because we are dealing with a local PDE, and all the information that in the original problem \eqref{eq:flp} would be coming from outside the unit ball is encoded in the values of $u$ in $B_1$ for $y>0$. To study the regularity of the problem, we will focus on a solution to (\ref{eq:symugp}-\ref{eq:symineq}).

Notice that the result of \cite{AC} is the particular case of this problem when $s=1/2$ and $\varphi=0$. The case $\varphi \neq 0$ does not seem to follow from the case $\varphi=0$ in a straight forward way. The main idea of writing the problem as a local-type partial differential equation is to be able to use usual techniques for regularity of partial differential equations, like monotonicity formulas and classification of blowup profiles.

The problem can be thought as the minimization of the functional
\[ J(v) = \int_{B_1} |\grad v|^2 \yed X \]
from all functions $v$ in the weighted Sobolev space $W^{1,2}(B_1,|y|^a)$ such that $v(x,0) \geq \varphi(0)$. Following the intuition explained in \cite{CS} this could be interpreted as an obstacle problem, where the obstacle is only defined in a set of codimension $1+a$, and $a$ is not necessarily an integer number.

The theory of degenerate elliptic equations in weighted Sobolev spaces plays an important role in this work, specially the results in the paper \cite{FKS}.

The strategy of our proof is the following. From \cite{S}, it is enough to show the regularity around free boundary points, so we assume that the origin $(0,0)$ belongs to the free boundary. We will show that $\int_{S_r} u^2 \yed \sigma$ decays in the appropriate way. To see that, we will study the limit as $r \to 0$ of a variation of the Almgren's formula:
\[ \Phi(r) \approx  \frac{r \int_{B_r} \abs{\grad u}^2 \yed X}{\int_{S_r} \abs{u}^2 \yed \sigma} \]

In order to obtain the possible limits of Almgren's formula, we make a blowup and it turns out that the limits correspond to the degrees of the global homogeneous solutions of (\ref{eq:symrn1}-\ref{eq:symrn4}) with $\varphi=0$. Finally. The minimum possible of such degrees of homogeneity can be found either using a monotonicity formula similar to the one from \cite{AC}, or using the sharp result of \cite{S} for the case of convex contact sets.

Later, the regularity of the free boundary is addressed using that in the blowup solutions the contact set is a half space.

\section{Preliminaries}
We will start with some elementary properties of the equation $\Dey{v}=0$. In the case $a=0$ most of these properties are very classical results for harmonic functions. Intuitively, these results generalize to the case $a \in (-1,1)$ through the idea of the fractional-dimension extension \cite{CS}. On a first reading, it may be a good idea to skim through this section quickly.

Since $|y|^a$ is an $A_2$ weight, the following proposition is a particular case of a result in \cite{FKS} that we will use in this paper
\begin{prop} \label{prop:fks}
Assume $\Dey{v} = 0 \text{  in } B_r$. Then there is an $\alpha>0$ such that the function $v$ is $C^\alpha(B_{r/2})$ and 
\begin{equation*}
\norm{v}_{C^\alpha(B_{r/2})} \leq \frac{C}{r^\alpha} \osc_{B_r} v
\end{equation*}
\end{prop}

For a proof of the above proposition check \cite{FKS}.

The Harnack inequality is also available.
\begin{prop} \label{prop:harnack}
Assume $\Dey{v} = 0 \text{  in } B_r$, then
\begin{equation*}
\sup_{B_{r/2}} u \leq C \inf_{B_{r/2}} u
\end{equation*}
Moreover, if instead we had a right hand side: $\Dey{v} = f(X) \text{  in } B_r$, we would have
\begin{equation} \label{eq:hrhs}
\sup_{B_{r/2}} u \leq C \inf_{B_{r/2}} u + C r^2 \sup_{B_r} |f|
\end{equation}
\end{prop}

The estimate with the right hand side \eqref{eq:hrhs} can be deduced from the weak-Harnack inequalities applied to $u + c|X|^2$ and $u + c(r^2 - |X|^2)$, as it is standard. Notice that this construction can be translated in the $x$ direction but not in $y$. We are goint to apply the estimate \eqref{eq:hrhs} only in balls centered on $\{y=0\}$ (in the proof of Corollary \ref{cor:nondeg}).

Using Proposition \ref{prop:fks} and the translation invariance of the equation in the $x$ variable, we obtain the following result.
\begin{prop} \label{prop:gradx}
Assume $\Dey{v} = 0 \text{ in } B_r(X_0)$ for some $r>0$. Then
\begin{align*}
\sup_{B_{r/2}(X_0)} |\grad_x v| &\leq \frac{C}{r} \osc_{B_r(X_0)} v \\
[\grad_x v]_{C^\alpha(B_{r/2}(X_0))} &\leq \frac{C}{r^{1+\alpha}} \osc_{B_r(X_0)} v \qquad \text{for some small $\alpha>0$}
\end{align*}
where write $[f]_{C^\alpha(D)}$ to denote the seminorm
\[ [f]_{C^\alpha(D)} = \sup_{x,y \in D} \frac{|f(x)-f(y)|}{|x-y|^\alpha} \]
\end{prop}

\begin{remark}
Notice that the estimate in Proposition \ref{prop:gradx} refers to the derivatives with respect to $x$, not with respect to $y$.
\end{remark}

\begin{proof}
Assume first $r=1$. From Proposition \ref{prop:fks} we know apriori that $v$ is $C^\alpha$ and
\[ \norm{v}_{C^\alpha(B_{3/4}(X_0))} \leq C \osc_{B_1(X_0)} v \]

Given a tangential unit vector $\tau = (\tau',0)$ The incremental quotient:
\[ v_h = \frac{v(X+h\tau ) - v(X)}{|h|^\alpha} \]
is bounded in $B_{7/8}$ independently of $h$. Moreover, by the linearity and translation invariance of the equation, $v_h$ is also a solution of $\Dey{v_h}=0$ and then by Proposition \ref{prop:fks}, it is $C^\alpha$.

This implies that $v$ is of class $C^{2\alpha}$ in the $x$ variable (see \cite{CC}, or \cite{MS} to see an example where the method is only applied in the tangential direction). Iterating the same argument a finite number of times we obtain that $v$ is Lipschitz and finally $C^{1,\alpha}$ in the tangential directions in $B_{1/2}$

The estimate of the Proposition follows by scaling.
\end{proof}

\begin{cor} \label{cor:dx}
Assume $\Dey{v} = 0 \text{ in } B_r(X_0)$ for some $r>0$. Then
\begin{align*}
\sup_{B_{r/2}(X_0)} |D^k_x v| &\leq \frac{C}{r^k} \osc_{B_r(X_0)} v \\
[D^k_x v]_{C^\alpha(B_{r/2}(X_0))} &\leq \frac{C}{r^{k+\alpha}} \osc_{B_r(X_0)} v \qquad \text{for some small $\alpha>0$}
\end{align*}
\end{cor}

\begin{proof}
Assume first $r=1$, the general case follows by scaling. By Proposition \ref{prop:gradx}, 
\begin{equation*}
\sup_{B_{1/2}(X_0)} |\grad_x v| \leq C \osc_{B_1(X_0)} v
\end{equation*}

Moreover, $w_j=\partial_{x_j} v$ is also a solution to  $\Dey{w_j} = 0$ in $B_1(X_0)$ since the equation is linear and translation invariant in $x$. Applying Proposition \ref{prop:gradx} again we obtain
\begin{align*}
\sup_{B_{1/4}(X_0)} |D_x^2 v| &= \sup_{B_{1/4}(X_0)} \sum_j |\grad_x w_j|\\
&\leq C \osc_{B_{1/2}(X_0)} w\\
&\leq C \osc_{B_1(X_0)} v
\end{align*}

For $k=2$, the proposition follows by scaling.

Iterating Proposition \ref{prop:gradx} further, we obtain the bounds for larger values of $k$ with the same reasoning.
\end{proof}

\begin{prop} \label{prop:dyy}
Assume $\Dey{v} = 0  \text{ in } B_r(X_0)$. Then for any $r\leq1$, 
\begin{equation*}
\sup_{B_{r/2}(X_0)} \abs{v_{yy} + \frac{a}{y} v_y} \leq \frac{C}{r^2} \osc_{B_r(X_0)} v
\end{equation*}
\end{prop}

\begin{proof}
From Corollary \ref{cor:dx} we have that
\[ \sup_{B_{r/2}(X_0)} |D^2_x v| \leq \frac{C}{r^2} \osc_{B_r(X_0)} v \]

However, from the equation we have that
\[ \lap_x v = - v_{yy} - \frac{a}{r} v_y \]

Therefore 
\[ \sup_{B_{r/2}(X_0)} |v_{yy} + \frac{a}{y} v_y| \leq \frac{C}{r^2} \osc_{B_r(X_0)} v \]
\end{proof}

The following is a Liouville type result

\begin{lemma} \label{lem:liouville}
Let $v$ be a global solution of
\[ \Dey{v(X)} = 0 \qquad X \in \R^n \times \R \]
such that $v(x,y) = v(x,-y)$ and $|v(X)| \leq C |X|^k$. Then $v$ is a polynomial.
\end{lemma}

\begin{proof}
We use Proposition \ref{prop:gradx} and Proposition \ref{prop:dyy} and induction in the degree $k$. The following elementary fact will be used: if $\grad_x v$ is a polynomial and $v(0,y)$ is a polynomial in $y$, then $v$ is a polynomial.

Let us start with the case $k \leq 1$. By taking $r \to +\infty$ in Proposition \ref{prop:dyy} we obtain $v_{yy} + \frac{a}{y} v_y = 0$. This gives a simple second order ODE for each $x$ whose solutions have the general form $b y |y|^{-a} + c$. Therefore, for a fixed $x$,  $v(x,y) = b y |y|^{-a} + c$, but since we assume that $v$ is symmetric in $y$, we have that $v(x,y)$ must be a constant for every fixed $x$.

On the other hand, taking $r \to +\infty$ in Corollary \ref{cor:dx} (with $k=2$) we have that $D_x^2 v = 0$. Therefore for each fixed $y$, $v(x,y)$ is a first order polynomial in $x$.

Combining the two facts above we have that if $k \leq 1$ then $v$ is a polynomial of the form $v(x,y) = bx + c$.

Now we consider larger values of $k$.

From Proposition \ref{prop:gradx}, $|\grad_x v(X)| \leq C |X|^{k-1}$. Moreover, $\grad_x v$ is also a global solution of the same equation symmetric in $y$. By the inductive hypothesis $\grad_x v$ is a polynomial of degree $k-1$.

From Proposition \ref{prop:dyy}, $|v_{yy} + \frac{a}{y} v_y| \leq C |X|^{k-2}$. Observe that
\[ v_{yy} + \frac{a}{y} v_y = |y|^{-a} \partial_y \left( |y|^a \partial_y v \right) \]
Thus $v_{yy} + \frac{a}{y} v_y$ satisfies the same equation as $v$ ($|y|^a v_y$ satisfies the conjugate equation. See \cite{CS}). Therefore $v_{yy} + \frac{a}{y} v_y$ is a polynomial of degree $k-2$.

In particular, for $x=0$, $v_{yy}(0,y) + \frac{a}{y} v_y(0,y)$ is some even polynomial $p(y)$ of degree at most $k-2$. Let us say that $p(y) = a_0 + a_2 y^2 + \dots + a_{2d} y^{2d}$. Then $v(0,y)$ must be
\[ v(0,y) = c + b y |y|^{-a} + \frac{a_0}{2(1+a)} y^2 + \frac{a_2}{4(3+a)} y^4 + \dots + \frac{a_2d}{(2d+2)(2d+1+a)} y^{2d+2} \]
but since $v$ is even in $y$, $b=0$ and $v(0,y)$ is a polynomial.

Since $v(0,y)$ is a polynomial in $y$ and $\grad_x v$ is a polynomial, then $v$ is a polynomial.
\end{proof}

\begin{remark}
The symmetry condition $v(x,y) = v(x,-y)$ is necessary. The simplest counterexample without that condition is $v(x,y) =|y|^{-a} y$.
\end{remark}

The equation $\Dey u=0$ can be understood as the Laplace equation in $n+1+a$ dimensions \cite{CS}. In that context, the following lemma is just the mean value theorem.

%Improve to the symmetrized one. It would make it simpler.
\begin{lemma} \label{lem:meanvalue}
Let $v$ be a function for which
\begin{equation*}
\Dey {v(X)} \leq C |y|^a \abs{X}^k  \text{  for } X \in B_1
\end{equation*}
Then \[ v(0) \geq \frac{1}{\omega_{n+a} r^{n+a}} \int_{S_r} v(X) \yed \sigma - C r^{k+2} \] for any $r < 1$. where
\[ \omega_{n+a} := \int_{S_1} \yed \sigma \]
\end{lemma}

\begin{proof}
The proof is essentially the same as the mean value property for harmonic functions. Let us consider the case $C=0$ first.

We consider the test function
\[\Gamma(X) = \max \left( \frac{C_{n,a}}{|X|^{n+a-1}} - \frac{C_{n,a}}{r^{n+a-1}} , 0 \right) \]
where $C_{n,a} = (n+a-1)^{-1} \omega_{n+a}^{-1}$. Note that $\Gamma$ is supported in $B_r$, $\Gamma \geq 0$ and $\Dey{\Gamma} (X) = -\delta_0 + \mu$, where $\mu$ is the measure supported on $S_r$ given by $r^{-n-a} w_{n+a}^{-1} |y|^a \dd \sigma$ (See \cite{CS}).
Thus
\begin{align*}
0 &\leq \int_{B_1} -\Dey{v(X)} \Gamma(X) \dd X \\
&= \int_{B_1} \grad v(X) \cdot \grad \Gamma(X) \yed X \\
&= \int_{B_1} - v(X) \Dey \Gamma (X) \dd X \\
&= v(0) - \frac{1}{w_{n+a} r^{n+a}} \int_{S_r} v(X) \yed \sigma
\end{align*}
which proves the Lemma when $C=0$. In the case $C>0$, we apply the above computation to $v(X) - C |X|^2 / (2(n+a+1))$.
\end{proof}

%Do not symmetrize
%\begin{prop} \label{prop:inhomog}
%If $g:\R^n \times [0,\infty) \to \R$ is a bounded function, the inhomogeneous equation
%\begin{equation*}
%\begin{aligned}
%w(x,0) &= 0 && \text{for } x \in \R^n\\
%y^{-a} \Dy{w(x,y)} &= g(x,y) && \text{for } y >0  
%\end{aligned}
%\end{equation*}
%has a solution $w$ such that $\lim_{y \to 0} y^a w_y(x,y)$ is a $C^\alpha$ function in $\R^n$ for any $\alpha < 1+a$.
%\end{prop}

%\begin{proof}
%\textbf{TO DO.} (I hope my guess $\alpha=1+a$ is correct. This will be used to estimate the regularity of an auxiliary function constructed to locallize the result of my thesis.)
%\end{proof}

We will also need Poincar\'e type inequalities in the context of weighted Sobolev spaces. The following is a classical Poincar\'e inequality whose proof can be found in \cite{FKS}.

\begin{lemma} \label{lemma:poincare}
For any function $v \in W^{1,2}(B_1,|y|^a)$ the following inequality holds
\[ \int_{S_r} |v(X) - \overline v|^2 \yed \sigma \leq C r \int_{B_r} \abs{\grad v(X)}^2 \yed X \]
where $\overline{v} = \frac{1}{\omega_{n+a} r^{n+a}} \int_{S_r} v(X) \yed \sigma$ and $C$ is a constant that depends only on $a$ and dimension.
\end{lemma}

\begin{remark}
In \cite{FKS}, the inequality is done with $\overline v$ being the average in the whole ball instead of the sphere. Indeed, the modification in the proof is straight forward given that a function $v$ in the weighted Sobolev space $W^{1,2}(B_1,|y|^a)$ has a trace in $L^2(S_1,|y|^a)$.
\end{remark}

Another form of Poincar\'e inequality that will come handy is the following.
\begin{lemma} \label{lem:yetanotherpoincare}
For any $r<1$ there is a constant $C>0$ (depending only on $r$, $a$ and dimension) such that given any function $v \in W^{1,2}(B_1,|y|^a)$ the following inequality holds
\[ \int_{S_1} |v(X) - v(rX)|^2 \yed \sigma \leq C r \int_{B_1} \abs{\grad v(X)}^2 \yed X \]
\end{lemma}

The proof is standard. It can be done for example integrating $\grad v$ along the lines $sX$ with $s \in (r,1)$, or using the compactness of the trace operators from $W^{1,2}(B_1,|y|^a)$ to $L^2(S_r,|y|^a)$.

Since the main difficulties in our problem appear in the behaviour of the solution near the free boundary, we will assume that the origin belongs to it.

Almgren's monotonicity formula was proved for the problem (\ref{eq:symugp}-\ref{eq:symineq}) in \cite{CS} if $\varphi=0$. We cannot reduce the problem to that case. Instead, assuming $\varphi \in C^{2,1}$, we let $\tilde u(x,y) = u(x,y) - \varphi(x) + \frac{\lap \varphi(0)}{2 (1+a)} y^2$ so that $\Dey{\tilde u} = 0$ at the origin. Denote by $\Lambda = \{\tilde u = 0\} = \{ u = \varphi\}$ the contact set. The function $\tilde u$ has the following properties:
\begin{align}
\tilde u(x,0) &\geq 0 && \text{for } x \in B_1^* \label{eq:obs}\\
\tilde u(x,y) &= \tilde u(x,-y) \label{eq:sym}\\
\Dey{\tilde u(x,y)} &= |y|^a (\lap \varphi(x) - \lap \varphi(0)) =: |y|^a g(x) && \text{for } (x,y) \in B_1 \setminus \Lambda  \label{eq:rhs}\\
\Dey{\tilde u(x,y)} &\leq  |y|^a g(x) && \text{for } (x,y) \in B_1  \label{eq:rhsineq}
\end{align}
where we define the function $g(x) = (\lap \varphi(x) - \lap \varphi(0))$ that is Lipschitz as long as $\varphi \in C^{2,1}$.

The only problem we have is that the right hand side in \eqref{eq:rhs} is not zero. However, $|y|^a g(x) = |y|^a (\lap \varphi(x) - \lap \varphi(0))$ that decays to $0$ as $x \to 0$. In fact
\begin{equation} \label{eq:righthandside}
|\Dey{\tilde u(x,y)}| \leq C |y|^a |x|
\end{equation}

We would expect a small variation of Almgren's monotonicity formula to work, or at least to remain bounded as $r \to 0$ when we apply it to $\tilde u$. The functions $u$ and $\tilde u$ have the same regularity respect to $x$, so we will prove our estimates in terms of $\tilde u$. In order to simplify the notation we will write $u$ from now on, but we mean $\tilde u$.

%\section{Elementary properties of the equations \eqref{eq:obs}-\eqref{eq:supersol}}
\begin{lemma} \label{lem:secondpoincare}
Let $u$ be a function for which \eqref{eq:obs}, \eqref{eq:sym}, \eqref{eq:rhs}, \eqref{eq:rhsineq} hold and $u(0)=0$. Then
\[ \int_{S_r} |u(X)|^2 \yed \sigma \leq C r \int_{B_r} \abs{\grad u(X)}^2 \yed X + C r^{6+a +n} \]
for a constant $C$ depending only on $a$, $n$ and $\norm{ \varphi }_{C^{2,1} }$.
\end{lemma}

\begin{remark}
The exponent in $r^{6+a+n}$ is a consequence of the chosen regularity for $\varphi$, in this case $C^{2,1}$.
\end{remark}

\begin{proof}
By Lemma \ref{lemma:poincare}, we have
\begin{equation} \label{eq:a0}
 \int_{S_r} |u(X)|^2 \yed \sigma \leq C r \int_{B_r} \abs{\grad u(X)}^2 \yed X + \overline{u} \int_{S_r} u(X)  \yed \sigma
\end{equation}

Therefore, to prove the lemma we need to find a suitable upper bound for
\[ \int_{S_r}  |u(X)| \yed \sigma \]

From Lemma \ref{lem:meanvalue} we have that
\begin{equation*}
 0=u(0) \geq \frac{1}{\omega_{n+a} r^{n+a}} \int_{S_r}  u(X)  \yed \sigma  -  C r^{3}
\end{equation*}
therefore
\begin{equation}  \label{eq:a1}
 \int_{S_r}  u^+(X)  \yed \sigma \leq \int_{S_r}  u^-(X)  \yed \sigma + C r^{3+ a +n} 
\end{equation}
Now we estimate $u^-(x,y)$ by integrating along the straight line $(x,0)$ to $(x,y)$ and applying Cauchy-Schwarz. Since $u(x,0) \geq 0$, 
\begin{align*}
 u^-(x,y) &= u^-(x,y) - u^-(x,0) \leq \int_0^y \abs{\grad u(x,t)} \dd t \\
&\leq \left( \int_0^y \abs{\grad u(x,t)}^2 t^a \dd t \right)^{1/2} \left( \int_0^y t^{-a} \dd t \right)^{1/2} \\
&\leq C y^{\frac{1-a}{2}} \left( \int_0^y \abs{\grad u(x,t)}^2 t^a \dd t \right)^{1/2}
\end{align*}

Now we integrate the above inequality on $S_r$ and apply Cauchy-Schwarz again. Notice that on $S_r$ we have $\dd \sigma = \frac{r}{|y|} dx$.
\begin{align*}
\int_{S_r} u^-(x,y) \yed \sigma &=  \int_{S_r} u^-(x,y) |y|^{a-1} r \dx \\
&\leq Cr \int_{S_r} |y|^{\frac{a-1}{2}}  \left( \int_0^y \abs{\grad u(x,t)}^2 t^a \dd t \right)^{1/2} \dx \\
&\leq Cr \left( \int_{B_r} \abs{\grad u(x,y)}^2 \yed X \right)^{1/2} \left(\int_{S_r} |y|^{a-1} \dx\right)^{1/2} \\
&\leq  Cr^{(n+1+a)/2} \left( \int_{B_r} \abs{\grad u(x,y)}^2 \yed X \right)^{1/2}
\end{align*}

Combining with \eqref{eq:a1} we obtain
\[ \int_{S_r} |u(X)| \yed \sigma \leq Cr^{(n+1+a)/2} \left( \int_{B_r} \abs{\grad u(x,y)}^2 \yed X \right)^{1/2} + C r^{3+ a +n} \]

Putting this estimate back in \eqref{eq:a0} we have
\begin{align*}
\int_{S_r} |u(X)|^2 \yed \sigma &\leq C r \int_{B_r} \abs{\grad u(X)}^2 \yed X + \frac{1}{w_{n+a} r^{n+a}} \left( \int_{S_r} |u(X)| \yed \sigma \right)^2 \\
&\leq C r \int_{B_r} \abs{\grad u(X)}^2 \yed X + C r^{6 + a +n}
\end{align*}
\end{proof}

\begin{cor} \label{cor:secondpoincare}
Let $u$ be a function for which \eqref{eq:obs}, \eqref{eq:sym}, \eqref{eq:rhs}, \eqref{eq:rhsineq} hold and $u(0)=0$. Then
\[ \int_{B_r} |u(X)|^2 \yed X \leq C r^2 \int_{B_r} \abs{\grad u(X)}^2 \yed X + C r^{7+a +n} \]
for a constant $C$ depending only on $a$, $n$ and $\norm{ \varphi }_{C^{2,1} }$..
\end{cor}

\begin{proof}
Take Lemma \ref{lem:secondpoincare} and integrate in $r$.
\end{proof}

%The proof of following energy type inequality can also be found in \cite{FKS}.
%\begin{lemma} %\label{lem:energy}
%Let $v$ be a nonnegative function for which
%\begin{align}
%\Dey {v(X)} &\geq C_0 |y|  \text{ for } X \in B_1\\
%%\fld v(x,y) &\geq 0 && \text{for } x \in B_1^*
%\end{align}
%Then 
%\[ \int_{B_{1/2}} \abs{\grad u}^2 \yed X \leq C_1 \int_{B_1} |u|^2 \yed X + C_2 C_0^2\]
%for constants $C_1$ and $C_2$ depending only on dimension and $a$.
%\end{lemma}

\section{Frequency formula}
A crucial ingredient for our blowup analysis is the monotonicity of a frequency formula of Almgren type for functions $u$ that solve (\ref{eq:obs}-\ref{eq:rhsineq}), and $u(0)=0$. In the special case that the right hand side in \eqref{eq:rhs} is zero, then the frequency formula takes the simple classical form as shown in \cite{CS}:
\[ \Phi(r) = \frac{r \int_{B_r} \abs{\grad u}^2 \yed X}{\int_{S_r} \abs{u}^2 \yed \sigma} \]

In order to account for the right hand side in \eqref{eq:rhs}, we have to modify the formula with a suitable lower order term. We define the function:
\begin{equation} \label{eq:functionF}
 F_u(r) = \int_{S_r} u(X)^2 \yed \sigma
\end{equation}

Notice that in terms of $F$ we have
\[ r \frac{\dd}{\dd r} \log(F_u(r))  = \frac{r \int_{B_r} \abs{\grad u}^2 \yed X}{\int_{S_r} \abs{u}^2 \yed \sigma} + n + a. \]
Thus, the classical frequency formula, shown in \cite{CS}, can be rephrased as that $r \frac{\dd}{\dd r} \log(F_u(r))$ is monotone nondecreasing in $r$. We will use the following modification
\begin{equation} \label{eq:frequency}
 \Phi_u(r) = (r+C_0 r^2 ) \frac{\dd}{\dd r} \log \max(F_u(r) , r^{n+a+4})
\end{equation}

As long as we are talking about only one function, we will write $\Phi = \Phi_u$. Whenever there might be some ambiguity we will use the subindices.

\begin{thm} \label{thm:frequencyformula}
For some small $r_0$ and a constant $C_0$ large enough (depending on $a$, $n$ and $\norm{\varphi}_{C^{2,1}}$),  The function $\Phi$ (given in \ref{eq:frequency}) is monotone nondecreasing for $r<r_0$.
\end{thm}

This result is an extension of the corresponding one in \cite{CS}. The proof is somewhat technical and we believe it is better to skip it on a first reading of this paper. In order not to distract the reader we postpone the proof to the appendix, so that we can quickly concentrate on the regularity of the obstacle problem in the following section.

\section{Local $C^{1,\alpha}$ estimates}

In \cite{S}, it is proved that the solutions of \eqref{eq:flp} are in general $C^{1,\alpha}$ for any $\alpha<s$. We can take this result to the context of (\ref{eq:symugp}-\ref{eq:symineq}) locallizing the problem using a cutoff function.

\begin{lemma} \label{lem:c1ainx}
Given a $C^2$ function $\varphi :B_1^* \to \R$, we consider a function $u : B_1 \to \R$ satisfying equations  (\ref{eq:obs}-\ref{eq:rhsineq}). Assume
\begin{align*}
|u(X)| &\leq C && \text{for } X \in B_1\\
%u_{ee}(x,0) &\geq -C && \text{for } x \in B_1^* \text{ and $e$ is any unit vector such that $e \cdot e_{n+1} = 0$}\\
|\varphi(x)| + |D^2 \varphi (x)| &\leq C && \text{for } x \in B_1^*
\end{align*}
Then we have the following estimates depending only on $C$, $a$ and $n$:
\begin{itemize}
\item $u_{ee}(x,0)$ is bounded below for $x \in B_{1/2}^*$ for any unit vector $e$ such that $e \cdot e_{n+1} = 0$.
\item $\lim_{y \to 0} |y|^a u_y(x,y) \in C^{\alpha} (B^*_{1/2})$ for any $\alpha < 1-s$. %is bounded in $B^*_{1/2}$.
\item $u(x,0) \in C^{1,\alpha} (B^*_{1/2})$ for any $\alpha < s$.
\end{itemize}
\end{lemma}

\begin{proof}
We apply here the results in \cite{S}. The only technical difficulty in order to apply immediately the main theorem therein is that our function $u$ has a bounded domain instead of being defined in the whole space. The solution is to use a smooth radially symmetric cutoff function $\eta : \R^{n+1} \to \R$ that vanishes outside $B_1$ and equals $1$ in $B_{2/3}$.

We want to show that $\tilde u(x) = \eta(x,0) u(x,0)$ satisfies the hypothesis of Proposition 5.7 in \cite{S}.

Since the cutoff function $\eta$ is radially symmetric, $\eta u$ keeps the same Neumann condition on the boundary $\R^n \times \{0\}$: 
\begin{align*}
\fld {(\eta u)} (x,y) &= 0 && \text{for } u(x,0)>0 \\
\fld {(\eta u)} (x,y) &\leq 0 && \text{for } x \in \R^n
\end{align*}
However, the function $\eta u$ is not necessarily a solution of $\Dey{\eta u} = 0$ in the half space $\R^n \times (0,+\infty)$. Let $f(X) = \Dey{ \left( \eta(X) u(X) \right)}$. We solve the following ``correction'' equation:
\begin{align*}
w(x,0) &= 0 && \text{for } x \in \R^n \\
\Dey{w(x,y)} &= f(x,y) && \text{for } y>0
\end{align*}

We look at the right hand side $f(X) = \Dey{\left( \eta(X) u(X) \right)} = u \Dey{\eta} + y^a \grad u \grad \eta$.
Notice that $f(X) = 0$ for $X \in B_{1/2}$. From Corollary \ref{cor:dx}, we have that $w$ is $C^\infty(B_{1/2}^*)$, and thus also is $\lim_{y \to 0} |y|^a \partial_y w(x,y)$.

% Recall that $\Dey \eta = y^a \left( \lap_X \eta + \frac{a}{y} \eta_y \right)$. Then \[ f(X) = y^a \left( u \lap_X \eta + u \frac{a}{y} \eta_y + 2 \grad u \grad \eta \right) \]Since $\eta$ is radially symmetric $\eta_y(x,0) = 0$, and $\eta_y$ is smooth. Therefore, from the above formula, $y^{-a} f(X)$ is bounded, and from Proposition \ref{prop:inhomog}, $\fld w(x,y)$ is $C^\alpha$ in $x$ for any $\alpha<1+a$. Notice also that $f(X)=0$ for $X \in B_{1/2}$. This also implies that $\fld w(x,y)$ is $C^\infty(B_{1/2}^*)$. 

Now the function $\eta u-w$ is a solution to 
\[ \Dey (\eta u - w)(X) = 0 \]
in the whole half space. Its restriction $\tilde u(x)=\eta(x,0) u(x,0) - w(x,0)$ satisfies
\begin{align*}
\tilde u &\geq 0 &&\\
% \tilde u_{ee} &\geq -c I && \text{for any unit vector } e \\
(-\lap)^s \tilde u(x) &\geq \phi(x) := -c \fld w(x,y) \\
(-\lap)^s \tilde u(x) &= \phi(x) && \text{for those $x$ where }  \tilde u(x) > 0
\end{align*}

Moreover, for $\varphi = (-\lap)^{-s} \phi$, $\tilde u+\varphi$ is a solution of \eqref{eq:flp}. In \cite{S} it is shown that $(-\lap)^s \tilde u$ is bounded and $\tilde u_{ee}$ is bounded below for any unit vector $e$. From the first one we obtain the boundedness of $\lim_{y \to 0} |y|^a u_y(x,y)$ and from the second one the bound below for $u_{ee}(x,0)$.

Applying Proposition 5.7 in \cite{S}, we get $\tilde u \in C^{1,\alpha}$ for any $\alpha < s$ in the $x$ variable and obtain the other estimate. Recall that $(-\lap)^s \tilde u(x) = \fld{\tilde u(x,y)}$, so the $C^\alpha$ estimate for $\lim_{y \to 0} |y|^a u_y(x,y)$ comes from the corresponding estimate for $(-\lap)^s \tilde u(x)$ in $\cite{S}$.
\end{proof}

\begin{remark}
Notice that the function $w$ that we construct in the proof above complies with all the hypothesis (4.1-4.3) and (4.32 - 4.34) from \cite{S} that are needed in order to apply Proposition 5.7 from there. In \cite{S} only obstacles $\varphi$ with compact support are considered for simplicity in the proof of existence of a solution, but the results extend whenever the problem makes sense.
\end{remark}

A $C^{1,\alpha}$ regularity estimate with respect to the $x$ variable only implies the same regularity estimates in both $x$ and $y$ for the case $a=0$. In any case, what we can say is that $y^a \partial_y u(x,y)$ remains uniformly bounded as $y \to 0$.

\begin{prop} \label{prop:c1ainboth}
Let $u$ be as in Lemma \ref{lem:c1ainx}. Then we have the following estimates
\begin{align}
\grad_x u(X) \in C^\alpha(B_{1/2}) && \text{ for any } \alpha < s \label{eq:lipinx} \\
|y|^a \partial_y u(X) \in C^\alpha(B_{1/2}) && \text{ for any } \alpha < 1-s \label{eq:whateveriny}
\end{align}
and the estimates on the corresponding norms depend only on the constants in the statement of Lemma \ref{lem:c1ainx}.
\end{prop}

\begin{proof}
From Lemma \ref{lem:c1ainx} we have both estimates at the boundary $\{y=0\}$. Since $w_i = \partial_{x_i} u$ also solves $\Dey w = 0$ in $B_1 \setminus \{y=0\}$, the estimate extends to the interior and we obtain \eqref{eq:lipinx}. On the other hand $w(x,y) = |y|^a \partial_y u(x,y)$ solves the equation $\dv (|y|^{-a} \grad w(X))=0$ (See \cite{CS}) and we obtain \eqref{eq:whateveriny}.
\end{proof}

The above proposition will be used together with the following compactness lemma to prove the existence of a blowup limit.

\begin{lemma} \label{lemma:compactness}
Let $v_j:B_1 \to \R$ be a bounded sequence of functions in $W^{1,2}(B_1,|y|^a)$. Assume there is a constant $C$ so that
\begin{align*}
\abs{\grad_x v_j(X)} &\leq C && \text{for } X \in B_1 \\
\abs{\partial_y v_j(X)} &\leq C |y|^{-a} && \text{for } X \in B_1
\end{align*}
and also that for each $\delta>0$ the sequence is uniformly $C^{1,\alpha}$ in $B_{1-\delta} \cap \{|y| \geq \delta\}$.

Then there is a subsequence $v_{j_k}$ that converges strongly in $W^{1,2}(B_1,|y|^a)$.
\end{lemma}

\begin{proof}
Since the sequence $v_j$ is bounded in $W^{1,2}(B_1,|y|^a)$, then there is a subsequence that converges strongly in $L^2(B_1,|y|^a)$ (see \cite{K}). 

Since for each $\delta>0$ the sequence is uniformly $C^{1,\alpha}$ in $B_{1-\delta} \cap \{|y| \geq \delta\}$, we can extract a subsequence so that $\grad v_j$ converges uniformly in $B_{1-\delta} \cap \{|y| \geq \delta\}$ for any $\delta>0$. Thus, $\grad v_j$ converges pointwise in $B_1 \setminus \{y=0\}$.

In order to get strong convergence in $W^{1,2}(B_1,|y|^a)$, we will show that each partial $\partial_k v_j$ converges in $L^2(B_1,|y|^a)$. Notice that each partial $\partial_k v_j$ already converges almost everywhere

We use two facts. The first is that if a sequence of functions $w_j$ is bounded in $L^p$ (for $1<p<\infty$) and converges almost everywhere, then it converges weakly in $L^p$. The second is that if $w_j$ is an almost everywhere and weakly convergent sequence of functions bounded in $L^p(\Omega)$ for some measure space $(\Omega,\mu)$ such that $\mu(\Omega) < +\infty$ and $p>2$ then $w_j$ converges strongly in $L^2(\Omega)$ (See \cite{MT}).

Let us do it first for $\partial_y v_j$. Since $|y|^a \abs{\partial_y v_j}$ is bounded and converges in $B_1 \setminus \{y=0\}$, then $|y|^a \partial_y v_j$ converges strongly in $L^2(B_1)$, and therefore $v_j$ converges strongly in $L^2(B_1,|y|^a)$.

With respect to $\grad_x v_j(X)$, we know that it is bounded uniformly in $L^\infty(B_1,|y|^a)$ and converges pointwise in $B_1 \setminus \{y=0\}$, therefore it converges strongly in $L^2(B_1,|y|^a)$.
\end{proof}

\begin{remark} \label{rk:noc1a}
Notice that the function $u$ that solves (\ref{eq:symugp}-\ref{eq:symineq}) can only be $C^{1,\alpha}$ in both variables $x$,$y$ in the case $a \leq 0$. If $a > 0$, then from the fact that $y^a \partial_y u(x,y)$ has a nonzero limit as $y \to 0^+$ for some points $x$ in the contact set $\{u=\varphi\}$, then $\partial_y u$ can never be bounded. However, after the change of variables $z = y |y|^{-a}$ (see \cite{CS}), it is to be expected that the function $u$ would be $C^{1,\alpha}$ in $x$ and $z$ for a certain value of $\alpha$.
\end{remark}

\section{Blowup profiles}

The purpose of this section is to characterize all posible blowup profiles of solutions to the equation. The problem we consider is to find all functions $u : \R^{n+1} \to \R$ with the properties:
\begin{align}
\text{$u$ is }&\text{homogeneous of degree $k$} \label{eq:ge0}\\
u(x,0) &\geq 0 && \text{for } x \in \R^n \label{eq:ge1} \\
u(x,y) &= u(x,-y) \\
\Dey {u(X)} &= 0 && \text{in } \R^n \times \R \setminus \{(x,0) : u(x,0)=0\}  \label{eq:ge2}  \\
\Dey{u(X)} &\leq 0 && \text{in } \R^n \times \R \text{ in the distributional sense.}  \label{eq:ge3}  \\
u_{\tau \tau} (X) &\geq 0 && \text{for any vector $\tau$ in $\R^n \times \{0\}$} \label{eq:ge4} 
\end{align}

The following proposition gives a lower bound for the degree of homogeneity $k$, and it is all that is needed for the optimal regularity result.

\begin{prop} \label{prop:globalprofiles}
If there is a solution $u$ of (\ref{eq:ge0}-\ref{eq:ge4}), then $k \geq (3-a)/2 = 1+s$, where $(1-a)/2 = s$
\end{prop}

\begin{proof}
The idea of this proof is to use the main theorem from \cite{S}. The theorem cannot be applied immediately because our function $u$ is not bounded, and even worse, in general its growth at infinity is too high to even define the fractional Laplacian as a distribution. The solution is to use a smooth radially symmetric cutoff function $\eta : \R^{n+1} \to \R$ like
\[ 
\eta(X) = \left\{ \begin{aligned}
1 & \text{ if } |X| < 1/2 \\
e^{-|X|^2} & \text{ if } |X|>1
                 \end{aligned} \right. \]

Now we want to show that $\tilde u(x) = \eta(x,0) u(x,0)$ satisfy the hypothesis of theorem 5.2 in \cite{S}. Notice that since $u$ is convex in the $x$ direction, then the contact set $\{(x,0) : \tilde u(x,0) = 0 \} = \{(x,0) : u(x,0)=0\}$ is convex. The proof follows as in Lemma \ref{lem:c1ainx} but applying Theorem 5.2 in \cite{S}. Applying that theorem we conclude that $u \in C^{1,s}$ in the $x$ variable, and therefore $k \geq 1+s$.
\end{proof}

\begin{remark}
An alternative approach to show Proposition \ref{prop:globalprofiles} could be to use a monotonicity formula similar to how it is done in \cite{AC}.
\end{remark}

Even though the proposition above is enough to obtain the optimal regularity of the solutions $u$ of (\ref{eq:symugp}-\ref{eq:symineq}), we need to fully characterize all solutions of (\ref{eq:ge0}-\ref{eq:ge4}) in order to study the regularity of the free boundary.

As before, we have $\Lambda = \{ (x,0) \in \R^n : u(x,0)=0 \}$. Let us also define
\begin{equation*}
\Lambda_* = \{ (x,0) \in \R^n : \fld u(x,y) < 0 \}
\end{equation*}
Notice that $\overline \Lambda_*$ is the support of $\Dey u$, and by \eqref{eq:ge2} $\Lambda_* \subset \Lambda$. Moreover, by the continuity of $u$, $\Lambda$ is closed and  $\overline{\Lambda_*} \subset \Lambda$.

\begin{lemma} \label{lem:polynomialprofile}
If $\Lambda_*$ has $H^n$-measure zero then $u$ is a polynomial of degree $k$.
\end{lemma}

\begin{proof}
From Proposition \ref{prop:c1ainboth} $|y|^a u_y(x,y)$ is locally bounded. If $\Lambda_*$ has $H^n$-measure zero, then \[\lim_{y \to 0} |y|^a u_y(x,y) = 0 \text{ a.e. in } x.\] Therefore for $\lim_{y \to 0} |y|^a u_y(x,y) = 0$ weak-$*$ in $L^\infty$. In \cite{CS}, it is shown that this implies that $u$ is a global solution of 
\[ \Dey{u} = 0 \qquad \text{in } \R^n \times \R \]

We conlude the proof using Lemma \ref{lem:liouville}.
\end{proof}

\begin{prop} \label{prop:1+s}
If $\Lambda$ has positive $H^n$-measure then either $u \equiv 0$, or $k = 1+s$ and $\Lambda$ is half of $\R^n$.
\end{prop}

\begin{proof}
We first observe that if $\Lambda_*$ has $H^n$-measure zero then $u \equiv 0$. Otherwise from Lemma \ref{lem:polynomialprofile} $u(x,0)$ would be a polynomial vanishing in a set of positive measure in $\R^n$, thus it would be constant zero. The polynomial $u$ must have the form
\[ u(x,y) = p_1(x) y^2 + \dots p_j(x) y^{2j} \]

Computing $\Dey{u}$ in terms of the above expresion, a simple iterative computation shows that $\Dey{u} = 0$ only if $p_1 = p_2 =\dots =p_j=0$.

Let us now consider the case when $\Lambda_*$ has a positive $H^n$-measure.

Since $u$ is homogeneous, $\Lambda_*$ is a cone. Let us assume that $e_n$ is a direction inside this cone. Since $\Lambda_*$ has a positive $H^n$-measure, then it must be a \emph{thick} convex cone in the sense that a neighborhood of $e_n$ is contained in $\Lambda_*$. Therefore, for any $x \in \R^n$, $x+he_n \in \Lambda_*$ for $h$ large enough.

Since 
\[ c \lim_{y \to 0} \frac{u(x,y) - u(x,0)} { |y|^{-a} y} = \lim_{y \to 0} |y|^a u_y(x,y) \leq 0 \]
then $u$ is not positive in a neighborhood of $e_n$. Therefore, for any $X \in \R^{n+1}$, $u(X+h e_n) \leq 0$ for $h$ large enough.

By hypothesis \eqref{eq:ge4}, $u$ is convex in the $e_n$ direction. The function $w = -u_{x_n}$ is decreasing and cannot be negative at any point $X$ because otherwise $\lim_{h \to +\infty} u(X+h e_n) = +\infty$ contradicting the above.

On the other hand, $w=0$ on $\Lambda$ and 
\[ \Dey{w(X)} = 0 \qquad \text{in } \R^{n+1} \setminus \overline{\Lambda_*} \supset \R^{n+1} \setminus \overline{\Lambda} \]

Thus, $w$ must be the first eighenfunction corresponding to minimizing the spherical integral
\[ \int_{S_r} |\grad_\theta w|^2 |y|^a \dd \theta \] 
from all functions $w$ such that $w = 0$ on $\Lambda$ and 
\[ \int_{S_r} |w|^2 |y|^a \dd \theta = 1 \] 

Since $\Lambda$ is convex, it covers at most half of the sphere $S_r \cap \{y=0\}$. If it was exactly half of the sphere then we have the explicit expresion
\[ w(x,y) = c \left( \sqrt{x_n^2 + y^2} - x_n \right)^s \]
which is, up to a multiplicative factor, the only positive solution of $\Dey{w}=0$ which vanishes in $\{y=0 \wedge x_n \geq 0\}$ (the computation is somewhat lenghty, but it can be done quickly with a computer algebra system).

Notice that the above explicit function is not a solution accross $\{y=0 \wedge x_n \geq 0\}$. Therefore for any convex cone $\Lambda$ that is strictly contained in $\{y=0 \wedge x_n \geq 0\}$, there must be another function which gives a smaller eigenvalue, and then a smaller degree of homogeneity $k$ than $1+s$. But from Proposition \ref{prop:globalprofiles}, $k \geq 1+s$. Therefore the only possibility is $k=1+s$. Moreover, we also see in the explicit function $w$ that $\Lambda$ has to be $\{y=0 \wedge x_n \geq 0\}$, which is half of $\R^n$. Recall that $e_n$ is an arbitrary direction inside the cone $\Lambda$.
\end{proof}

\begin{prop} \label{prop:uniqueblowup}
Up to rotations and multiplicative constants, there is a unique solution of (\ref{eq:ge0}-\ref{eq:ge4}) that is homogeneous of degree $1+s$. 

For this solution the free boundary is flat, there is a unit vector $e$ such that $\Lambda = \{ (x,0) : x \cdot e \geq 0 \}$, and $\partial_e u = c \left( \sqrt{(x \cdot e)^2 + y^2} - (x \cdot e) \right)^s$.
\end{prop}

\begin{proof}
In the proof of Proposition \ref{prop:1+s}, once $u_{x_n}$ is uniquely determined as
\[ w = -u_{x_n} = c \left( \sqrt{x_n^2 + y^2} - x_n \right)^s \]
we have $\Lambda = \{ u = 0\} = \{y=0 \wedge x_n \geq 0\}$, and integrating on the lines parallel to $e_n$ we determine $u(x,0)$ for every $x$. Now, if we had two solutions $u_1$ and $u_2$ homogeneous of degree $1+s$, coinciding on $\{y=0\}$, then necessarily $u_1(x,y)-u_2(x,y) = c |y|^s y$ for some constant $c$ and $y>0$. But that constant must be zero in order for $u_1$ and $u_2$ to be solutions accross $\{y=0\} \setminus \Lambda$.

Replacing $e_n$ by an arbitrary unit vector $e$ normal to $\bdary \Lambda$ we get the rest of the result.
\end{proof}

\begin{cor} \label{cor:twovariables}
Let $u$ be the solution of (\ref{eq:ge0}-\ref{eq:ge4}) that is homogeneous of degree $1+s$ and such that $e_n$ is normal to its free boundary, then $u$ is constant with respect to $x_1, x_2, \dots, x_{n-1}$.
\end{cor}

\begin{proof}
Translating $u$ in any direction $\tau$ orthogonal to $x_n$ and $y$ would give another global profile with the same free boundary. By uniqueness, $u$ must be invariant in those directions.
\end{proof}

We can summarize the results of the section in the following Theorem
\begin{thm} \label{thm:globalprofiles}
Let $u$ be a solution of (\ref{eq:ge0}-\ref{eq:ge4}). There are two posibilities:
\begin{enumerate}
\item The degree $k = 1+s$. $\Lambda$ is a half space, and $u$ depends only on two directions.
\item The degree $k$ is an integer greater or equal to $2$, $u$ is a polynomial, and $\Lambda$ has $H^n$-measure zero.
\end{enumerate}
\end{thm}

\section{Blowup sequences and optimal regularity}
The optimal regularity of the solution will be obtained by carefully analyzing the possible values of $\Phi(0)$. Recall that when $\Phi$ is constant, its value depends on the degree of homogeneity of the function. In this case $\Phi(0)$ represents the assymptotic behaviour at the origin.
%Now we come back to a solution $v$ of \eqref{eq:obs}-\eqref{eq:supersol}.

\begin{lemma} \label{lem:Phi0}
Let $u$ be a solution to (\ref{eq:obs}-\ref{eq:rhsineq}). Then either $\Phi_u(0) = n+a + 2(1+s)$ or $\Phi_u(0) \geq n+a + 4$
\end{lemma}

In order to prove Lemma \ref{lem:Phi0} we will consider the following scaled version of $F$ that represents the growth of $u$ at the origin.
\[ d_r = \left( r^{-(n+a)} \int_{S_r} |u|^2 \yed X \right)^{1/2} = \left( r^{-(n+a)} F_u(r) \right)^{1/2} \]

Using the value of $d_r$ we define the following nonhomogenous blowup sequence that will also be useful for studying the regularity of the free boundary in the next sections.
\begin{equation} \label{eq:blowupsequence}
u_r(X) = \frac{1}{d_r} u( r X)
\end{equation}

\begin{lemma} \label{lem:blowup} Let $u$ be as in Lemma \ref{lem:Phi0}.
 If $\liminf_{r \to 0} \frac{d_r}{r^2} = +\infty$, then there is a sequence $r_k \to 0$ and a nonzero function $u_0 :\R^{n+1} \to \R$ such that
\begin{align}
 u_{r_k} &\to u_0 && \text{in } H^1(B_{1/2}) \\
 u_{r_k} &\to u_0 && \text{uniformly in } B_{1/2} \\
 \grad_x u_{r_k} &\to \grad_x u_0 && \text{uniformly in } B_{1/2} \\
 |y|^a \partial_y u_{r_k} &\to |y|^a \partial_y u_0 && \text{uniformly in } B_{1/2}
\end{align}

Moreover, $u_0$ satisfies the conditions (\ref{eq:ge0}-\ref{eq:ge4}) and its degree of homogeneity is $(\Phi_u(0) - n - a)/2$.
\end{lemma}

\begin{proof}
The function $u_r$ is constructed in a way so that $\norm{u_r}_{L^2(S_1, |y|^a)} = 1$ for every $r$. Since $\liminf_{r \to 0} \frac{d_r}{r^2} = +\infty$, then $F_u(r) > r^{n+a+4}$ for $r$ small enough. Let us also consider the values of $r$ smaller than the $r_0$ of Theorem \ref{thm:frequencyformula}.

We will first show that $u_r$ remains bounded in $W^{1,2}(B_1,|y|^a)$ using the monotonicity of the frequency formula. We have that
\begin{equation} \label{eq:o0}
\begin{aligned}
 \Phi(r_0) &\geq \Phi(r) \\
&\geq (r + C_0 r^2) \frac{\dd}{\dd r} \log \max( F_u(r), r^{n+a+4} ) \\
&\geq  (r + C_0 r^2)\frac{2 \int_{S_r} u u_\nu \yed \sigma}{\int_{S_r} |u|^2 \yed \sigma} + (n+a)(1+C_0r) \qquad \text{(since $F_u(r) > r^{n+a+4}$)}\\
&\geq 2(r + C_0 r^2)\frac{\int_{B_r}  u \Dey{u} + |\grad u|^2 \yed X}{\int_{S_r} |u|^2 \yed \sigma} + (n+a)(1+C_0r)\\
\end{aligned}
\end{equation}

Recalling that $|\Dey{u}| \leq C |y|^a |x|$ in $\R^{n+1} \setminus (\Lambda \times \{0\})$, and that $u=0$ on $\Lambda \times \{0\}$, we estimate
\begin{equation*} %\label{eq:o1}
 \abs{\int_{B_r} u  \Dey{u} \dd X} \leq C r^{\frac{n+a+3}{2}} \left( \int_{B_r} |u|^2 \yed X \right)^{1/2}
\end{equation*}

By Corollary \ref{cor:secondpoincare}, 
\begin{equation*} 
 \int_{B_r} |u(X)|^2 \yed X \leq C r^2 \int_{B_r} \abs{\grad u(X)}^2 \yed X + C r^{7+a +n} 
\end{equation*}
Thus
\begin{equation} \label{eq:o1}
 \abs{\int_{B_r}  u \Dey{u} \dd X} \leq C r^{\frac{n+a+5}{2}} \left( \int_{B_r} |\grad u|^2 \yed X \right)^{1/2} + C r^{n+a+5} 
\end{equation}

From Lemma \ref{lem:secondpoincare}, 
\[  \int_{B_r} \abs{\grad u(X)}^2 \yed X \geq \frac{1}{r} \int_{S_r} |u(X)|^2 \yed \sigma - C r^{5+a +n} \]
but since we have $d_r / r^2 \to +\infty$,
\[ \int_{S_{r}} |u(X)|^2 \yed \sigma \geq \tilde C r^{n+a+4} \]
for a constant $\tilde C$ as large as we wish as $r \searrow 0$.  Thus
\begin{align}
  \int_{B_{r}} \abs{\grad u(X)}^2 \yed X &\geq \tilde C r^{n+a+3} - C r^{5+a +n} \\
&\geq \frac{\tilde C}{2} r^{n+a+3} \qquad \text{for $r$ small } \label{eq:o2}
\end{align}

Putting it all back in \eqref{eq:o1},
\[ \abs{\int_{B_{r}}  u \Dey{u} \dd X} \leq C r^{n+a+4} \qquad \text{for $r$ small} \]

Comparing with \eqref{eq:o2} we see that for $r$ small
\[ \int_{B_r}  u \Dey{u} + |\grad u|^2 \yed X \geq \frac{1}{2} \int_{B_r} |\grad u|^2 \yed X \]

Now we continue with \eqref{eq:o0},
\[ \begin{aligned}
\Phi(r_0) &\geq (r + C_0 r^2)\frac{\int_{B_{r}} |\grad u|^2 \yed X}{\int_{S_{r}} |u|^2 \yed \sigma} + (n+a)(1+C_0r)\\
&\geq \frac{1}{2} r \frac{\int_{B_{r}} |\grad u|^2 \yed X}{\int_{S_{r}} |u|^2 \yed \sigma} + \frac{(n+a)}{2} \qquad \text{for $r$ small}
&\geq \frac{1}{2} \int_{B_1} |\grad u_r|^2 \yed X
   \end{aligned}
\]

Combining the above inequality with the fact that $\int_{S_1} |u_r|^2 \yed \sigma = 1$ and Poincar\'e inequality, we obtain that the sequence $u_r$ remains bounded in $W^{1,2}(B_1,|y|^a)$.

Since $u_r^+=\max(u_r,0)$ and $u_r^-=\max(-u_r,0)$ are subsolutions of the equation
\[ \dv( |y|^a \grad u) \geq -C r |y|^a |x| \]
then we have that $u_r$ is bounded in $L^\infty(B_{3/4})$ (See \cite{FKS}).

We also know from Lemma \ref{lem:c1ainx} that $u$ is semiconvex in $x$, or in other words that $\partial_{ee} u \geq -C$ for any tangential unit vector $e$. This implies a bound of the same type for all the sequence $u_r$ since we have the choice of scaling $\lim_{r \to 0} \frac{d_r}{r^2} = \infty$.

Notice that the functions $u_r$ are solutions of a uniformly elliptic equation with smooth coefficients in $B_1 \cap \{y>\delta/2\}$ for any $\delta > 0$. Then clearly $u_r$ are uniformly $C^{1,\alpha}$ in $B_{1-\delta} \cap \{y>\delta\}$ for any $\alpha,\delta>0$. We can therefore apply Proposition \ref{prop:c1ainboth} and Lemma \ref{lemma:compactness} to obtain a subsequence $u_{r_j}$ such that $u_{r_j}$ converges strongly in $W^{1,2}(B_{1/2},|y|^a)$ to some function $u_0$.

From Theorem \ref{thm:frequencyformula}, we have that $\Phi_u(r)$ is monotone and converges to some value as $r \to 0$.  Moreover, the computations right above show that as $r \to 0$
\begin{equation} \label{a1}
\begin{aligned}
\Phi_u(rs) &\approx rs \frac{\int_{B_{rs}} |\grad u|^2 \yed X}{\int_{S_{rs}} |u|^2 \yed \sigma} + (n+a) \\
&\approx r \frac{\int_{B_{r}} |\grad u_s|^2 \yed X}{\int_{S_{r}} |u_s|^2 \yed \sigma} + (n+a)
\end{aligned} 
\end{equation}

Now we let $s=r_j \to 0$, to obtain
\begin{equation} \label{eq:almgrenLimit}
 r \frac{\int_{B_{r}} |\grad u_0|^2 \yed X}{\int_{S_{r}} |u_0|^2 \yed \sigma} = \Phi_u(0) - (n+a)
\end{equation}

In order to pass to the limit in the above expression we need $u_s \to u_0$ in $W^{1,2}(B_r,|y|^a)$ (which we have), and we also need that the denominator remains bounded away from zero (i.e. that $u_0 \neq 0$). When $r > 1-\delta$, if $\delta$ is small, this is a consequence of Lemma \ref{lem:yetanotherpoincare}. In other words, we have that if $\delta$ is small enough
\[ \int_{S_{(1-\delta)r}} |u_s|^2 \yed \sigma \geq c \int_{S_{r}} |u_s|^2 \yed \sigma \]
for a constant $c$ depending on $\delta$, $\Phi(1)$ and dimension.

We can iterate this inequality $k$ times so that $(1-\delta)^k < r$ and obtain a uniform bound from below for the denominator \eqref{a1}, so that we can pass to a nonzero limit and get \eqref{eq:almgrenLimit}.

Since $\lim_{r \to 0} \frac{d_r}{r^2} = +\infty$, we have that for any unit vector $\tau = (\tau',0)$
\[ \partial_{\tau \tau} u_r = \frac{r^2}{d_r} u_{\tau \tau} ( r X) \geq -C \frac{r^2}{d_r} \to 0 \]
So, in the limit, $u_0$ is convex in the $x$ direction.

Each function $u_r$ is a solution to
\begin{align}
u_r(x,0) &\geq 0 && \text{for } x \in B_1^* \label{eq:ur0} \\
\Dey{u_r(X)} &= \frac{r^{2-a}}{d_r} \Dey{u} (rX) = \frac{r^2}{d_r} |y|^a g(rx) && \text{for } X \in B_1 \setminus \{(x,0) : u_r(x,0) = 0 \} \label{eq:a9} \\
\Dey{u_r(X)} &\leq \frac{r^2}{d_r} |y|^a g(x)  && \text{for } X \in B_1 \label{eq:ur2} 
\end{align}
Notice that the right hand side in \eqref{eq:a9} goes to $0$ as $r \to 0$ since
\[ \abs{ \frac{r^2}{d_r} |y|^a g(rx) }\leq C \frac{r^2}{d_r} |y|^a |r x| \to 0 \]

Therefore $u_0$ is a solution of the homogeneous problem
\begin{align*}
u_0(x,0) &\geq 0 && \text{for } x \in B_1^* \\
\Dey{u_0(X)} &= 0  && \text{for } X \in B_1 \setminus \{(x,0) : u_0(x,0) = 0 \} \\
\Dey{u_0(X)} &\leq 0  && \text{for } X \in B_1
\end{align*}

For this problem, Almgren's frequency formula applies in the usual way without an error correction \cite{CS}. So, from \eqref{eq:almgrenLimit}, we conclude that $u_0$ is homogeneous in $B_{1/2}$, and its degree of homogeneity is exactly $(\Phi_u(0)-(n+a))/2$. Since it is homogeneous, then it can be extended to $\R^{n+1}$ as a global solution of the homogeneous problem. 

The aditional fact that we can choose $r_k$ so that $u_{r_k}$, $\grad_x u_{r_k}$ and $|y|^a \partial_y u_{r_k}$ converge uniformly comes from the apriori estimates of Proposition \ref{prop:c1ainboth}.
\end{proof}

\begin{proof}[Proof of Lemma \ref{lem:Phi0}]
Two things may happen:
\[ \liminf_{r \to 0} \frac{d_r}{r^2} \left\{ \begin{array}{l}
= +\infty  \qquad \text {first case} \\
< +\infty \qquad \text {second case}
\end{array} \right.\]

For the \textbf{first case}, we use Lemma \ref{lem:blowup} to find the blowup profile $u_0$. Then we apply Theorem \ref{thm:globalprofiles} to obtain that the degree of homogeneity of $u_0$ is $1+s$ or at least $2$, and thus $\Phi_u(0) = \Phi_{u_0} (0) = n + a +2 (1+s)$ or $\Phi_u(0) = \Phi_{u_0} (0) \geq n+a+4$.

Now we turn to the \textbf{second case}.

If $F_u(r_j) < r_j^{n+a+4}$ for a sequence $r_j \to 0$, then $\Phi_u(r_j)$ equals to $n+a+4$ for those values of $r_j$ and $\Phi_u(0) = n+a+4$.

On the other hand, let us assume that $F_u(r) \geq r^{n+a+4}$ for $r$ small. Since we are considering the \emph{second case}, we have
\[ r_j^{n+a+4} \leq F_u(r_j) \leq C r_j^{n+a+4} \]
for some constant $C$ and a sequence $r_j \searrow 0$. Taking logs in the above inequality
\[ (n+a+4) \log r_j \leq \log F_u(r_j) \leq C + (n+a+4) \log r_j \]
We will show that $\Phi_u(0) = \lim_{j \to +\infty} \Phi_u(r_j) \geq n+a+4$ in this case. Suppose otherwise that for small $r_j$, $\Phi_u(r_j) \leq n+a+4-\eps_0$, then we take $r_m<r_n<<1$ and obtain
\begin{align*}
(n+a+4) (\log r_n - \log r_m) - C&\leq \log F_u(r_n) - \log F_u(r_m)\\
&\leq \int_{r_m}^{r_n} \frac{\dd}{\dd r} \log F_u(r) \dd r \\
&\leq (n+a+4-\eps_0) (\log r_n - \log r_m)
\end{align*}
which gives a contradiction when we take $\log r_n - \log r_m \to +\infty$. This finishes the proof.
\end{proof}

\begin{prop} \label{prop:blowup} Let $u$ be as in Lemma \ref{lem:Phi0}.
Assume $\Phi(0)=n+a+2(1+s)$. There is a family of rotations $A_r$ such that $u_r \circ A_r$ converges to the unique global profile $u_0$ of degree $1+s$ (as mentioned in Proposition \ref{prop:uniqueblowup}). The convergence is in the sense that
\begin{align}
 u_r \circ A_r &\to u_0 && \text{in } W^{1,2}(B_{1/2},|y|^a) \\
 u_r \circ A_r &\to u_0 && \text{uniformly in } B_{1/2} \\
 \grad_x (u_r \circ A_r) &\to \grad_x u_0 && \text{uniformly in } B_{1/2} \\
 |y|^a \partial_y (u_r \circ A_r) &\to |y|^a \partial_y u_0 && \text{uniformly in } B_{1/2}
\end{align}
\end{prop}
\begin{remark}
 By a rotation $A$, naturally what we mean is an orthonormal matrix, i.e. $A^t A = Id$. The set of orthonormal matrices is compact, therefore Proposition \ref{prop:blowup} immediately implies the existence of convergent subsequences.
\end{remark}
\begin{proof}
In the proof of Lemma \ref{lem:Phi0}, we saw that in the case when $\liminf_{r \to 0} \frac{d_r}{r^2} < +\infty$ we have $\Phi(0) \geq n+a+4$. If $\Phi(0)=n+a+2(1+s)$ we must have $\liminf_{r \to 0} \frac{d_r}{r^2} = +\infty$.

Consider $\norm{v} = \norm{v}_{H^1(B_{1/2},|y|^a)} + \norm{v}_{C^0(B_{1/2})} + \norm{\grad_x v}_{C^0(B_{1/2})} + \norm{|y|^a \partial_y v}_{C^0(B_{1/2})}$.

From Lemma \ref{lem:blowup}, we can find a sequence $r_k \to 0$ such that $u_{r_k}$ converges to a global profile of degree $1+s$. Notice that $\int_{S_1} |u_0|^2 \yed X = 1$, and from Proposition \ref{prop:uniqueblowup}, $u_0$ is unique up to rotation. In particular
\[ \lim_{k \to \infty} \inf_{\text{all rotations } A } \norm{ u_0 - u_{r_k} \circ A } = 0 \]
where we consider $\norm{v} = \norm{v}_{H^1(B_{1/2},|y|^a)} + \norm{v}_{C^0(B_{1/2})} + \norm{\grad_x v}_{C^0(B_{1/2})} + \norm{|y|^a \partial_y v}_{C^0(B_{1/2})}$.

If there was any sequence $r_k$ for which $\norm{ u_0 - u_{r_k} \circ A }$ stays away from zero, then by Lemma \ref{lem:blowup}, we can find a subsequence $r_{k_j} \to 0$ such that $u_{r_{k_j}}$ converges to a global profile of degree $1+s$. Notice that $\int_{S_1} |u_0|^2 \yed \sigma = 1$, and from Proposition \ref{prop:uniqueblowup}, $u_0$ is unique up to rotation. In particular
\[ \lim_{j \to \infty} \inf_{\text{all rotations } A } \norm{ u_0 - u_{r_{k_j}} \circ A } = 0 \]
arriving to a contradiction. Therefore $\norm{ u_0 - u_{r_k} \circ A }$ cannot stay away from zero for any sequence $r_k$, and the result follows.
\end{proof}

\begin{lemma} \label{lem:l2toca}
Let $F$ be as in \eqref{eq:functionF}. If 
\begin{equation} \label{eq:a10}
F(r) \leq C r^{n+a+2(1+\alpha)}
\end{equation}
for any $r<1$, then $u(0) = \abs{\grad u(0)} = 0$, and $u$ is $C^{1,\alpha}$ at the origin in the sense that
\[ |u(X)| \leq C' |X|^{1+\alpha} \]
for $|X|<1/2$ and a constant $C'$ depending only on $C$, $a$ and dimension.
\end{lemma}

\begin{proof}
Consider $u^+ = \max (u,0)$ and $u^- = \max(-u,0)$. We have
\begin{align*}
\Dey{u^+} &\geq |y|^a g(x) \geq -C |y|^a |x|\\
\Dey{u^-} &\geq -|y|^a g(x) \geq -C |y|^a |x|
\end{align*}

For some $r>0$, let $u_0$ be the $\Dey{}$-harmonic replacement of $u^+$ in $B_r$. We note that 
\[ 0=\Dey{u_0} \leq \Dey{ \left(u^+ + C \frac{  |X|^2 - r^2 }{2(n+1+a)} \right) } \]
then by comparison principle \[ u_0 \geq u^+ - C' r^2 \]

By \eqref{eq:a10}, we have
\[ \int_{S_r} |u_0(X)|^2 \yed \sigma(X) = \int_{S_r} |u^+(X)|^2 \yed \sigma(X) \leq C r^{n+a+2(1+\alpha)} \]

Since $|y|^a$ is an $A_2$ weight, we use the estimates in \cite{FKS} to conclude
\[ \sup_{B_{r/2}} u_0 \leq C r^{1+\alpha} \]

Then $u^+(X) \leq u_0(X)+ r^2 \leq C r^{1+\alpha}$ for any $X \in B_{r/2}$. We can use the same reasoning for $u^-$ and conclude the proof.
\end{proof}

\begin{lemma}\label{lem:l2decay}
Let $\Phi_u$ be the frequency formula as in \eqref{eq:frequency}. Assume that $\Phi_u(r) \searrow \mu$ as $r \searrow 0$. Then
\[F_u(r) \leq C r^\mu\]
for any $r<1$ and a constant $C$ depending only on $F_u(1)$ and the constant $C_0$ in \eqref{eq:frequency}.
\end{lemma}

\begin{proof} 
For simplicity, we omit the subscripts and write $\Phi(r)$ and $F(r)$. 

Let $\tilde F(r) = \max(F(r),r^{n+a+4})$. By definition $\tilde F \geq F$. From the fact that $\Phi$ is monotone nondecreasing
\[ \mu = \Phi(0) \leq \Phi(r) = (r+C_0 r^2) \frac{\dd}{\dd r} \log \tilde F(r) \]

then
\[ \frac{\dd}{\dd r} \log \tilde F(r) \geq \frac{\mu}{r+C_0 r^2} \]
and integrating we obtain
\begin{align*}
\log \tilde F(1) - \log \tilde F(r) &\geq \int_r^1 \frac{\mu}{s+C_0 s^2} \dd s \\
&\geq -\mu \log r - \mu \int_r^1 \frac{1}{s} - \frac{1}{s+C_0 s^2} \dd s \\
&\geq -\mu \log r - \log(1-C_0) \mu
\end{align*}
and then $\log \tilde F(r) \leq \mu \log r + C_0 \mu + \log \tilde F(1)$. Exponentiating the two sides of the inequality we obtain
\[ F(r) \leq \tilde F(r) \leq C r^\mu \]
where $C = \tilde F(1) (1+C_0)^\mu$.
\end{proof}

\begin{thm} \label{th:decay}
Let $u$ be a function for which \eqref{eq:obs}, \eqref{eq:sym}, \eqref{eq:rhs}, \eqref{eq:rhsineq} hold and $u(0)=0$. Then
\[ u(X) \leq C |X|^{1+s} \sup_{B_1} |u| \]
for some constant $C$ depending only on $n$, $a$, and the Lipschitz norm of $g$.
\end{thm}

\begin{proof}
From Lemma \ref{lem:Phi0}, $\mu=\Phi(0) \geq n + a + 2(1+s)$, so the Theorem is proved by applying directly Lemma \ref{lem:l2decay} and Lemma \ref{lem:l2toca}.
\end{proof}

\begin{cor}
If $\varphi \in C^{2,1}$, the solution $u$ to the problem \eqref{eq:flp} is in the class $C^{1,s}$.
\end{cor}

\begin{proof}
Using the equivalence between the problem \eqref{eq:flp} and (\ref{eq:symrn1}-\ref{eq:symrn4}), Theorem \ref{th:decay} shows that $u - \varphi$ has a $C^{1,s}$ decay at free boundary points. As shown in \cite{S}, this is enough to prove that $u \in C^{1,s}$.
\end{proof}

\begin{cor}
If the function $u$ solves the problem (\ref{eq:obs}-\ref{eq:rhsineq}), then $u(x,y_0) \in C^{1,s}(B_{1/2}^*)$ for any value $0 < y_0 < 1/2$.
\end{cor}

\begin{cor}
If the function $u$ solves the problem (\ref{eq:symugp}-\ref{eq:symineq}), then $u(x,y_0) \in C^{1,s}(B_{1/2}^*)$ for any value $0 < y_0 < 1/2$.
\end{cor}

\begin{remark}
As in Remark \ref{rk:noc1a}, it is not true that $u$ is $C^{1,s}$ in the whole variable $X$ for general values of $a$. Interestingly enough however, the right decay in $y$ takes place at free boundary points. In the variables $x$ and $z$ (as in Remark \ref{rk:noc1a}) we would not have the right decay at free boundary points, and therefore the optimal regularity would not be uniform in both variables.
\end{remark}

\begin{remark}
The fact that the $C^{1,\alpha}$ estimates for $\alpha$ small happen naturally in the $x,z$ variables and the right sharp decay at the free boundary points happens in the $x,y$ variables may somehow reflect the respective nondivergent and divergent nature of each of these two results.
\end{remark}

\section{Free boundary regularity}

In this section we study the regularity of the free boundary. First we show it is a Lipschitz surface around nonsingular points. Then we will apply the boundary Harnack principle to tangential derivatives in order to show that it is a $C^{1,\alpha}$ surface. This is the same approach that was used in \cite{ACS} and also \cite{C1}.

\subsection{Lipschitz continuity of the free boundary}

\begin{thm} \label{thm:lipschitzfb}
Assume $\mu = \Phi(0) < n+a+4$. Then there exists a neighborhood of the origin $B_\rho$ and a tangential cone $\Gamma'(\theta, e_n) \subset \R^n \times \{0\}$ such that, for every $\tau \in \Gamma'(\theta,e_n)$ we have $D_\tau u \geq 0$. In particular, the free boundary $F(u)$ is the graph of a Lipschitz function $x_n = f(x_1,\dots, x_{n-1})$.
\end{thm}

The theorem will follow by applying the following lemma to a tangential derivative $h = D_\tau u_r$, where $u_r$ is the blow up family that defines the limiting global profiles, for $r$ small. % (also to keep small the constant $\gamma$ appearing in the lemma). Notice that by the properties of the limiting profile we can choose the opening $\theta$ of the cone as large as we like. This implies that the positive tangential derivatives have a nondegeneracy porperty of the type \[D_\tau u(X) \geq c d_X^\beta \] with some $\beta < 2$, where $d_X = \dist(X,F(u))$.

\begin{lemma} \label{lem:approximation}
Let $\Lambda$ be a subset of $\R^n \times \{0\}$. Assume $h$ is a continuous function with the following properties:
\begin{enumerate}
\item $|\Dey{h}| \leq \gamma |y|^a $
\item $h \geq 0$ for $|y| \geq \sigma > 0$, $h=0$ on $\Lambda$.
\item $h \geq c_0 > 0$ for $|y| \geq \frac{\sqrt{1+a}}{8n}$
\item $h \geq -\omega(\sigma)$ for $|y| < \sigma$, where $\omega$ is the modulus of continuity of $h$.
\end{enumerate}

There exists $\sigma_0 = \sigma_0(n,a,c_0,\omega)$ and $\gamma_0 = \gamma_0(n,a,c_0,\omega)$ such that, if $\sigma<\sigma_0$, $\gamma<\gamma_0$ then $h \geq 0$ in $B_{1/2}$.
\end{lemma}

\begin{proof}
Suppose $X_0 = (x_0,y_0) \in B_{1/2}$ and $h(X_0) = 0$. Let
\[ Q = \set{ (x,y) : |x-x_0| < \frac{1}{3}, |y| <\frac{\sqrt{a+1}}{4n} } \]
and \[ P(x,y) = |x-x_0|^2 - \frac{n}{a+1} y^2 .\]

Observe that $\Dey{P}=0$. Define \[ v(X) = h(X) + \delta P(X) - \frac{\gamma}{2(a+1)} y^2 .\]

Then:
\begin{align*}
v(X_0) &= h(X_0) + \delta P(X_0) - \frac{ \gamma }{ 2(a+1) } y_0^2 < 0 \\
\Dey{v} &= \Dey{h} + \delta \Dey{P} - \gamma |y|^a < 0 \leq 0 && \text{outside } \Lambda\\
v(X) &\geq 0 && \text{on } \Lambda
\end{align*}
Thus, $v$ must have a minimum on $\partial Q$. On
\[ \partial Q \cap \set{ |y| > \frac{ \sqrt{1+a}} {8n} } \]
we have \[v \geq c_0 - \frac{\delta}{16n} - \frac{\gamma}{32 n} > 0 \]
if $\delta$ and $\gamma$ are small depending on $c_0$ and $n$.

On \[ |x-x_0| = \frac{1}{3} \qquad \sigma \leq |y| \leq \frac{ \sqrt{1+a} }{8n} \]
one has \[ v \geq \delta \left( \frac{1}{9} - \frac{1}{64n} \right) - \frac{\gamma}{128 n^2} > 0 \]
if $\gamma$ is small compared to $\delta$.

Finally, on \[ |x - x_0| = \frac{1}{3} \qquad |y| < \sigma \] one has \[v \geq -\omega(\sigma) + \delta \left( \frac{1}{9} - \frac{n^2} {1+a} \sigma^2 \right) - \frac{\gamma \sigma^2}{2 (1+a)} > 0 \]
if $\sigma$ is small, depending on $\delta$ and $n$.

Hence, $v \geq 0$ on $\partial Q$ and we have a contradiction. Therefore $h \geq 0$ in $B_{1/2}$.
\end{proof}

Theorem \ref{thm:lipschitzfb} follows applying the Lemma \ref{lem:approximation} to a tangential derivative $h = D_\tau u_r$, where $u_r$ is the blow up family that defines the limiting global profiles, for $r$ small (also to keep small the constant $\gamma$ appearing in the lemma). This is a standard procedure by now that mimics the proof of Lipschitz regularity of the free boundary in the classical obstacle problem \cite{C1}. Below we include the proof for completeness.

Notice that by the properties of the limiting profile we can actually choose the opening $\theta$ of the cone as large as we like.

\begin{proof} [Proof of Theorem \ref{thm:lipschitzfb}]
We have $\mu = \Phi(0) < n+a+4$. By Lemma \ref{lem:Phi0}, $\mu=n+a+2(1+s)$. Moreover, the blowup sequence $u_r$ converges (up to subsequence) to the global profile $u_0$ of degree $1+s$ for which the free boundary is flat.

Let us assume that $e_n$ is the normal to the free boundary of $u_0$. We know by Proposition \ref{prop:uniqueblowup} that \[\partial_n u_0 = c \left( \sqrt{x_n^2 + y^2} - x_n \right)^s\]

For some $\theta_0>0$, let $v$ be any direction orthogonal to $y$ and $x_n$ such that $|v| < \theta_0$. From Theorem \ref{thm:globalprofiles} we have that $u_0$ is constant in the direction $v$, therefore if $\tau=e_n+v$, $\partial_\tau u_0 = \partial_n u_0$.

We know by Proposition \ref{prop:blowup} that $\grad_x u_r$ converges uniformly to $\grad_x u_0$ (up to subsequence). Therefore, for any $\delta_0$, there is an $r$ for which $|\partial_v u_\tau - \partial_\tau u_0| < \delta_0$ for any $\tau$ constructed as above ($\delta_0$ depends on $\theta_0$). We also have that $\partial_\tau u_0 = \partial_n u_0 = c \left( \sqrt{x_n^2 + y^2} - x_n \right)^s$. If we differentiate (\ref{eq:ur0}-\ref{eq:ur2}), we obtain the equations for $u_\tau$
\begin{equation} \label{eq:utau}
\Dey{(\partial_\tau u_r(X))} = \frac{r^2}{d_r} |y|^a r \partial_\tau g(rx) \leq C r |y|^a \qquad \text{for } X \in B_1 \setminus \{(x,0) : u_r(x,0) = 0 \},
\end{equation}
 and its right hand side goes to zero as $r \to 0$. Therefore for a small enough $r$, $\partial_\tau u_r$ will satisfy the hypothesis of Lemma \ref{lem:approximation}. Thus $\partial_\tau u_r$ will be nonnegative in $B_{1/2}$. This implies that for any point $x$ in the free boundary of $u$ such that $|x| < r/2$, then $x+\lambda \tau$ stays on one side of the free boundary for any positive value of $\lambda$, and thus the free boundary is Lipschitz.
\end{proof}

\subsection{Nondegeneracy properties}
Assuming as before that the origin is a free boundary point, we cannot assure any minimal growth of $u$ around $0$ without additional assumptions. It is easy to write down examples for the case $a=0$ in two dimensions. Even for homogeneous global solutions vanishing on half of $\R^n$ we have the following family of solutions of the thin obstacle problem expressed in polar cordinates.
\[ u = r^{2k+1} \sin( (2k+1) \theta ) \]
which grow arbitrarily slowly away from the origin.

The nondegenerate point that we will consider are those such that $\Phi(0)=(n+a)+2(1+s)$. At those points we know the exact asymptotic picture from Proposition \ref{prop:blowup}. For some small $r$, $u_r$ will be close to the blowup profile, and we will be able to apply the following lemmas.

\begin{lemma} \label{lem:nondeginy}
 There are positive numbers $\eps_0 = \eps_0(a,n)$ and $\delta_0 = (\frac{1}{12n})^{1/2s}$ so that the following is true:

Let $v$ be a function such that
\begin{align}
 \Dey{v(X)} &\leq \eps_0 && \text{for } X \in B_1^* \times (0,\delta_0)  \\
 v(X) &\geq 0 && \text{for } X \in B_1^* \times (0,\delta_0)  \\
 v(x,\delta_0) &\geq \frac{1}{4n} && \text{for } x \in B_1^*
\end{align}
Then $v(x,y) \geq C|y|^{2s}$ in $B_{1/2}^* \times [0,\delta_0]$.
\end{lemma}

\begin{proof}
This lemma refines the result of Lemma \ref{lem:approximation}, and therefore we can expect the proof to follow more or less the same pattern.

Let $X_0 = (x_0,y_0) \in B_{1/2}^* \times [0,\delta_0]$. We compare $v$ with the function
\[ w(x,y) = (1+\frac{\eps_0}{2}) y^2 - \frac{|x-x_0|^2}{n} + y^{2s} \]

We observe that inside the set $B_1^* \times (0,\delta_0)$, $\Dey{w} = \eps_0 \geq \Dey{v}$. We must check at the boundary now. For $y=0$, $w(x,0) \leq 0 \leq v(x,y)$. For $y=\delta_0$, \[w(x,\delta_0) \leq (1+\frac{\eps_0}{2}) \delta_0^2 + \delta_0^{2s} \leq 3 (\delta_0)^{2s} \leq \frac{1}{4n} \leq v(x,\delta_0).\]
And finally, where $|x-x_0| = 1/2$, $w(x,y) \leq 3 \delta_0^{2s} - \frac{1}{4n} \leq 0 \leq v(x,y)$. Therefore by comparison we have that $w \leq v$ in the whole set $B_1^* \times [0,\delta_0]$. In particular
\[ v(x_0,y_0) \geq w(x_0,y_0) \geq y^{2s} \]
\end{proof}

\begin{cor} \label{cor:nondeg}
Let $u$ be a solution of (\ref{eq:obs}-\ref{eq:rhsineq}). Moreover, let us assume that the function $g$ in \eqref{eq:rhs} is small in the sense that $|g(x)|,|\grad_x g(x)|\leq \eps_0$. Let $u_0$ be a global profile as in Proposition \ref{prop:uniqueblowup} such that $|\grad_x u -\grad_x u_0| \leq \eps_0$. Then if $\eps_0$ is small enough, there is a constant $c$ depending only on $n$ and $a$ such that we have the following lower bound for the tangential derivatives $u_\tau$
\[ u_\tau(X) \geq c \dist(X,\Lambda)^{2s} \]
for any $X \in B_{1/8}$ and $\tau$ being a tangential unit vector such that $|\tau - e_n| < 1/2$.
\end{cor}

\begin{proof}
From Lemma \ref{lem:approximation} and \eqref{eq:utau}, we know that $u_\tau$ will be positive in $B_{1/2}$.

Applying Lemma \ref{lem:nondeginy} to an appropiate multiple of $u_\tau$ we immediately get that
\[ u_\tau(x,y) \geq c |y|^{2s} \]
for every $(x,y) \in B_{1/4}$. Now, let $X = (x,y) \in B_{1/8}$ and $d = \dist(X,\Lambda)$. Consider the ball $B_{d/2}(X)$, then for the point $T=(x_T,y_T)$ at the top of this ball we have $y_T \geq d/2$. Therefore
\[ u_\tau(T) \geq c d^{2s} \]

By Harnack inequality (Proposition \ref{prop:harnack}) we know that $u_\tau(X)$ and $u_\tau(T)$ are comparable since $u_\tau$ is a solution of the equation in the ball $B_d(X)$. Therefore also
\[ u_\tau(X) \geq c d^{2s} \]
\end{proof}

\subsection{$C^{1,\alpha}$ regularity of the free boundary}
As in \cite{C1} or \cite{ACS}, we intend to prove a $C^{1,\alpha}$ regularity estimate for the free boundary applying the boundary Harnack inequality to two positive tangential derivatives $u_\tau$. We must be careful in order to handle the right hand side in (\ref{eq:obs}-\ref{eq:rhsineq}).

The first part in boundary Harnack (the Carleson estimate) can be proved even with a right hand side if we have a nondegeneracy condition like the one given by Corollary \ref{cor:nondeg}.

We use the notations from \cite{ACS}
\begin{itemize}
\item $D \subset B_1$ is a $\Dey{}-NTA$ domain.
\item $\Omega = \partial D \cap B_1$. If $Q \in \Omega$, $A_r(Q)$ is a point such that $B_{\eta r}(A(Q)) \subset B_r(Q) \cap D$. In particular $|A_r (Q) - Q| \sim r$.
\end{itemize}

Notice that the definition of $NTA$ domain given in \cite{ACS} differs slightly from the one used in \cite{FKJ}. A uniform capacity condition is assumed for the complement, instead of containing a full ball.

\begin{lemma} \label{lem:carleson}
Let $D \subset B_1$ be a $\Dey{}-NTA$ domain. Suppose $u$ is a positive function in $D$, continuously vanishing on $\Omega$. Assume that:
\begin{align}
|\Dey{u}|&\leq C|y|^a\\
u(X) &\geq c d_X^\beta && \text{for some $\beta$, } 0<\beta<2, \text{ where } d_X = \dist(X,\Omega) \qquad \text{(nondegeneracy)} \label{eq:snd}
\end{align}

The for every $Q \in \Omega \cap B_{1/2}$ and $r$ small:
\begin{equation}
\sup_{B_r(Q)} u \leq C(n,a,D) u(A_r(Q)) \label{eq:s2}
\end{equation}
\end{lemma}

\begin{proof}
Let $u^*$ be the $\Dey{}$-harmonic replacement of $u$ in $B_{2r}(Q) \cap D$, $r$ small. Then by thorem 6 in \cite{ACS} adapted to the operator $\Dey{}$ (the proof is identical using the results in \cite{FKS}) we get
\begin{equation} \label{eq:s1}
u^* (X) \leq C u^*(A_r(Q)) \qquad \text{in } B_r(Q) \cap D
\end{equation} 

On the other hand,
\begin{align*}
u^*(X) + C(|X-Q|^2 - r^2) &\leq u(X) && \text{on } \bdary( D \cap B_{2r} )\\
L_a{ \left( u^*(X) + C(|X-Q|^2 - r^2) \right) } &= C|y|^a \geq L_a u(X) && \text{in } D \cap B_{2r},
\end{align*}
then by comparison principle $u^* - u \leq r^2$ in $D \cap B_{2r}$. Similarly we obtain the other inequality and
$|u(X) - u^*(X)| \leq C r^2$.

From \eqref{eq:s1}
\[ u(X) \leq C( u(A_r (Q)) + Cr^2 ).\]

From the nondegeneracy property, $u(A_r(Q)) \geq cr^\beta$, and since $\beta < 2$, \eqref{eq:s2} follows.
\end{proof}

We will not prove the second part of Boundary Harnack in so much generality. Instead, we notice that for our particular case, we can obtain it from Lemma \ref{lem:approximation}.

\begin{lemma} \label{lem:bh}
 Let $D = B_1 - \Lambda$, where $\Lambda$ is a subset of $\R^n \times \{0\}$ such that there is a constant $\kappa$ so that the $L_a$-capacity of $B_r(X) \cap \Lambda$ is at least $\kappa r^{n+a-1}$ for every $X \in \Lambda$ and $r<1/2$ (\emph{this implies that $D$ is an NTA-domain as in \cite{ACS}}).

Let $u,v$ be two positive functions in $D$, as in Lemma \ref{lem:carleson}. Let us also assume that $u$ and $v$ are symmetric in $y$. Then there is a constant $c$ (depending on $n$, $a$, and the constant $\kappa$ above) such that
\begin{equation} \label{eq:bh}
 \frac{u(X)}{v(X)} \leq c \frac{u(e_{n+1})}{v(e_{n+1})} \qquad \text{in } B_{1/2}
\end{equation}
Moreover, the ratio $u/v$ is a $C^\alpha$ function in $B_{1/2}$ (uniformly up to the boundary $\Lambda$).
\end{lemma}

\begin{proof}
We divide both functions by their values at $e_{n+1}$, so that we can assume $u(e_{n+1}) = v(e_{n+1}) = 1$.

Using Lemma \ref{lem:carleson} and Harnack inequality, for any $\delta>0$ we have
\begin{align*}
 u(X) &\leq C && \text{for all $X \in B_{3/4}$}\\
 v(X) &\geq c && \text{for all $X \in B_{3/4} \cap \{y>\delta\}$}
\end{align*}
which means that for a constant $s>0$ small enough, $v-su$ satisfies the assumptions of Lemma \ref{lem:approximation}. Therefore $v-su$ is positive in $B_{1/2}$, or in other words $u(X)/v(X) \leq s$.

By a standard iterative argument \eqref{eq:bh} implies that $u/v$ is $C^\alpha$ up to the boundary for an appropriately small $\alpha$.
\end{proof}

The $C^{1,\alpha}$ regularity of the free boundary follows by applying Lemma \ref{lem:bh} to two positive tangential derivatives. We need the obstacle to be $C^{2,1}$.

\begin{thm} \label{thm:c1afb}
 Let $u$ be a solution of (\ref{eq:obs}-\ref{eq:rhsineq}). Assume $0$ is in its free boundary and $\Phi(0) < n+a+4$. Then the free boundary is a $C^{1,\alpha}$ $(n-1)$-dimensional surface around $0$.
\end{thm}

\begin{proof}
The proof follows by a standard argument in the subject \cite{C1}, \cite{CSa}. 

From Theorem \ref{thm:lipschitzfb}, the free boundary is Lipschitz. Moreover there is a cone of tangential directions $\tau$, so that $u_\tau$ is positive in a neighborhood of the origin. Let us assume that $e_n$ is an axis of that cone. We express the tangential derivatives of the level sets of $u$ as a quotient $u_\tau / u_{e_n}$ and use Lemma \ref{lem:bh} to show that all level sets up to the free boundary are uniformly $C^{1,\alpha}$ surfaces in the tangential direction.

We must show that we have all the necessary conditions to apply Lemma \ref{lem:bh}. First of all, we consider a rescaled function $u_r$ so that $u_r$ is close enough to the global profile $u_0$ and the right hand side in \eqref{eq:a9} is small. 

If $\tau$ is in the cone of directions from the free boundary where the directional derivatives are positive, then $u_\tau$ is positive. Otherwise we express
\[ \frac{u_\tau}{u_{e_n}} = 4 \frac{ u_{ e_n + \tau / 4 } }{u_{e_n}} - 4 \]
so we consider $e_n + \tau / 4$ instead, and reduce the problem for those $\tau$ where $u_\tau$ is positive.

From Corollary \ref{cor:nondeg}, we have the nondegeneracy condition for both $u_{e_n}$ and $u_\tau$. In order to apply Lemma \ref{lem:bh} we are only left to show that the contact set $\Lambda = \{(x,0): u(x,0)=0\}$ has a uniform $L_a$-capacity. But that is automatic since the free boundary is Lipschitz.
\end{proof}

\section*{Appendix: Frequency formula.}

In this appendix we prove Theorem \ref{thm:frequencyformula}. We recall the notation from \eqref{eq:functionF}
\begin{equation*}
 F_u(r) = \int_{S_r} u(X)^2 \yed \sigma
\end{equation*}
where $u$ is a function that solves  \eqref{eq:obs}, \eqref{eq:sym}, \eqref{eq:rhs}, \eqref{eq:rhsineq} , and $u(0)=0$.

As it was mentioned before, in the case of zero right hand side in \eqref{eq:rhs}, the frequency formula takes the simplest form and was proved in \cite{CS}. In general we will show that a constant $C_0$ exists so that
\begin{equation*}
 \Phi(r) = (r+C_0 r^2 ) \frac{\dd}{\dd r} \log \max(F_u(r) , r^{n+a+4})
\end{equation*}
is monotone increasing.

In order to prove that $\Phi$ is monotone, we will need the following lemma
\begin{lemma} \label{lem:almgreen}
The following identity holds for any $r \leq 1$.
\begin{equation} \label{eq:diva}
r \int_{S_r }  \left( \abs{u_\tau}^2 - \abs{ u_\nu }^2 \right) \yed \sigma = \int_{B_r} \left( (n+a-1) \abs{\grad u}^2 - 2 \langle X , \grad u \rangle g(x) \right) \yed X
\end{equation}
where $u_\tau$ is the gradient in the tangential direction to $S_r$ and $u_\nu$ is the derivative in the normal direction.
\end{lemma}

\begin{proof}
The Lemma is obtained by applying the divergence theorem to the vector field \[|y|^a \frac{\abs{ \grad u }^2 }{2} X - |y|^a \langle X , \grad u \rangle \grad u\]

We have to notice that
\[ \dv \left( y^a \frac{\abs{ \grad u }^2 }{2} X - |y|^a \langle X , \grad u \rangle \grad u \right) = 
y^a \frac{n+a-1}{2} \abs{ \grad u }^2  - \langle X , \grad u \rangle L_a u \]

Since $\langle X , \grad u \rangle$ is a continuous function on $B_r^*$ that vanishes on $\Lambda = \{u = 0\}$, we have that $\langle X , \grad u \rangle L_a u$ has no sigular part and coincides with $ \langle X , \grad u \rangle g(x)$.
\end{proof}

We now prove Theorem \ref{thm:frequencyformula}.

\begin{proof}[Proof of Theorem \ref{thm:frequencyformula}]
Notice that $\log F(r)$ is differentiable for $r>0$ since \[ \frac{\dd}{\dd r} \log F(r) = \frac{F'(r)}{F(r)} = \frac{2 \int_{S_r} u u_\nu \yed \sigma}{\int_{S_r} u^2 \yed \sigma} + \frac{n+a}{r} \]

On the other hand $\log r^{n+a+4}$ clearly has $\frac{n+a+4}{r}$ as derivative. When we take the maximum in \eqref{eq:frequency} it may happen that we get a nondifferentiable function, but in any case the jump in the derivative will be in the positive direction. In order to prove that $\Phi$ is monotone, we can concentrate in each of the two values for the maximum separately.

In the case $\Phi(r) = (r+C_0 r^2 ) \frac{\dd}{\dd r} \log  r^{n+a+4}$, we have $\Phi'(r) = (1+2 C_0 r) (n+a+4)$. Therefore $\Phi$ is monotone in this case.

When $\Phi(r) = (r+C_0 r^2 ) \frac{\dd}{\dd r} \log F(r)$ we need to work a little more. We will concentrate in this case now. Notice that this case happens only if $F(r) \geq  r^{n+a+4}$.

We expand the expression for $\Phi(r)$
\[ \Phi(r) = (r+C_0 r^2) \frac{ 2 \int_{S_r} u u_\nu \yed \sigma}{ \int_{S_r} u^2 \yed \sigma} + (n+a)(1+C_0 r) \]

The second term $ (n+a)(1+C_0 r)$ is clearly increasing, so it would suffice if we can show that the first term is nondecreasing. We compute its logarithmic derivate and show that it is nonnegative,
\begin{equation} \label{eq:aa}
\begin{aligned}
\frac{\dd}{\dd r} \log & \left( (r+C_0 r^2) \frac{ 2 \int_{S_r} u u_\nu \yed \sigma}{ \int_{S_r} u^2 \yed \sigma}  \right) =\\
&= \frac{\dd}{\dd r} \left( \log(r) + \log (1+C_0 r)  + \log  \int_{S_r} u u_\nu \yed \sigma - \log { \int_{S_r} u^2 \yed \sigma} \right) \\
%&= \frac{1}{r} + \frac{\eps r^{\eps-1}}{1+r^\eps} + \frac{\frac{\dd}{\dd r} \int_{S_r} u u_\nu \yed \sigma }{ \int_{S_r} u u_\nu \yed \sigma } - \frac{  \int_{S_r} 2 u u_\nu  \yed \sigma + \frac{n+a}{r}  \int_{S_r} u^2 \yed \sigma }  { \int_{S_r} u^2 \yed \sigma} \\
&= \frac{1}{r} + \frac{C_0}{1+C_0 r} + \frac{\frac{\dd}{\dd r} \int_{S_r} u u_\nu \yed \sigma }{ \int_{S_r} u u_\nu \yed \sigma } - \frac{  \int_{S_r} 2 u u_\nu  \yed \sigma}  { \int_{S_r} u^2 \yed \sigma} - \frac{n+a}{r}
\end{aligned}
\end{equation}

To estimate $\frac{\dd}{\dd r} \int_{S_r} u u_\nu \yed \sigma$ we notice that
\begin{equation} \label{eq:a2}
\begin{aligned}
 \int_{S_r} u u_\nu \yed \sigma  &= \int_{B_r} \dv ( |y|^a u \grad u ) \dd X \\
&= \int_{B_r} |y|^a |\grad u|^2 + u L_a u \dd X
%&= G(r) + \int_{B_r} u \dv ( y^a \grad u ) \dd X
\end{aligned}
\end{equation}
Thus
\begin{equation} \label{eq:a3}
\begin{aligned}
\frac{\dd}{\dd r} \int_{S_r} u u_\nu \yed \sigma &= \int_{S_r} |y|^a |\grad u|^2 + u \Dey{u}  \dd \sigma \\
&\geq \int_{S_r} |\grad u|^2 \yed \sigma -  C r^{\frac{n+a+2}{2}} \left( \int_{S_r} |u|^2 \yed \sigma \right)^{1/2} \\
&\geq \int_{S_r} |\grad u|^2 \yed \sigma -  C r^{\frac{n+a+2}{2}} F(r)^{1/2}
\end{aligned}
\end{equation}

Now we use Lemma \ref{lem:almgreen} to estimate $\int_{S_r} |\grad u|^2 \yed \sigma$.
\begin{align*}
 \int_{S_r} & |\grad u|^2 \yed \sigma = 2 \int_{S_r} |u_\nu|^2 \yed \sigma +  \frac{1}{r} \int_{B_r} \left( (n+a-1) \abs{\grad u}^2 - 2 \langle X , \grad u \rangle g(x) \right) \yed X \\
&= 2 \int_{S_r} |u_\nu|^2 \yed \sigma + \frac{n+a-1}{r}  \int_{S_r} u u_\nu \yed \sigma - \frac{1}{r} \int_{B_r} ( (n+a-1) u   +  2 \langle X , \grad u \rangle ) g(x) \yed X \\
\end{align*}

Recall that $|g(x)| \leq C |x|$. Putting the above estimate back in \eqref{eq:a3} we obtain
\begin{align*}
 \frac{\dd}{\dd r} \int_{S_r} u u_\nu \yed \sigma \geq  2 \int_{S_r} |u_\nu|^2 \yed \sigma &+ \frac{n+a-1}{r}  \int_{S_r} u u_\nu \yed \sigma \\
&- C r^{\frac{n+1+a}{2}} \left( \sqrt { G(r) }+ r \sqrt{ H(r)} + r^{\frac{1}{2}} \sqrt { F(r) } \right)
\end{align*}
where
\begin{align*}
F(r) = \int_{S_r} u(X)^2 \yed \sigma \\
G(r) = \int_{B_r} u(X)^2 \yed X \\
H(r) = \int_{B_r} |\grad u(X)|^2 \yed X
\end{align*}

Putting it all together in \eqref{eq:aa}, we write $\frac{\dd}{\dd r} \log \Phi(r) = R(r) + S(r)$, with
\[ R(r) = \frac{ 2 \int_{S_r} |u_\nu|^2 \yed \sigma }{ \int_{S_r} u u_\nu \yed \sigma } - \frac{  \int_{S_r} 2 u u_\nu  \yed \sigma}  { \int_{S_r} u^2 \yed \sigma} \geq 0 \] 
and
\[ \begin{aligned}
S(r) &= \frac{C_0}{1+C_0 r} - C r^{\frac{n+1+a}{2}} \frac{  \sqrt { G(r) }+ r \sqrt{ H(r)} + r^{\frac{1}{2}} \sqrt { F(r) } } {\int_{S_r} u u_\nu \yed \sigma} \\
&\geq \frac{C_0}{1+C_0 r} - C r^{\frac{n+1+a}{2}} \frac{  \sqrt { G(r) }+ r \sqrt{ H(r)} + r^{\frac{1}{2}} \sqrt { F(r) } } {H(r) -  r^{\frac{n+a+2}{2}} \sqrt{G(r)} }
\end{aligned} \]

In order to estimate $S$ properly, we must find bounds for $F$, $G$ and $H$. Recall that we are in the case when $F(r) \geq r^{n+a+4}$. From Lemma \ref{lem:secondpoincare} we know
\[ r^{n+a+4} \leq F(r) \leq C r H(r) + C r^{6+a+n} \]
and also, integrating the inequality in $r$ we have
\[ G(r) \leq C r^2 H(r) + C r^{7+a+n} \]

This means that, for small enough $r_0$, if $r<r_0$, 
\begin{align*}
F(r) &\leq C r H(r) \\
G(r) &\leq C r^2 H(r)
\end{align*}
and then
\[S(r) \geq \frac{C_0}{1+C_0 r} - C r^{\frac{n+1+a}{2}} \frac{ r \sqrt{ H(r)}} {H(r) -  r^{\frac{n+a+4}{2}} \sqrt{H(r)} }
\]

But since $F(r) \geq r^{n+a+4}$, also $H(r) \geq r^{n+a+3}$. We thus have, for $r < r_0$,
\[\begin{aligned}
S(r) &\geq \frac{C_0}{1+C_0 r} - C r^{\frac{n+1+a}{2}} \frac{ r } {\sqrt{H(r)} } \\
&\geq \frac{C_0}{1+C_0 r} - C r^{\frac{n+1+a}{2} + 1 - \frac{n+a+3}{2}} \\
&\geq \frac{C_0}{1+C_0 r} - C
\end{aligned}
\]

Which is positive if we choose $C_0$ large enough. This shows that $\Phi'(r) \geq 0$ for small $r$. If $\Phi$ is not monotone for larger values of $r$ in $(0,1)$, we just take $C_0$ larger enough so that $\Phi$ is nondecreasing for $0<r<1$.
\end{proof}

\bibliographystyle{plain}   % Here the bibliography 
\bibliography{css}             % is inserted.
\index{Bibliography@\emph{Bibliography}}%
\end{document}